\documentclass[11pt]{article}
\usepackage{hyperref}
\usepackage{color}
\usepackage{amssymb,amsmath, amsthm}
\newcommand{\bfi}{\bfseries\itshape}

\begin{document}
\newtheorem{thm}{Theorem}[section]
\newtheorem{lemma}{Lemma}[section]
\newtheorem{cor}{Corollary}[section]
\newtheorem{prop}{Proposition}[section]
\newtheorem{defn}{Definition}[section]
\newtheorem{rem}{Remark}[section]
\newtheorem{exam}{Example}[section]
%
\newcommand{\todo}[1]{\vspace{5 mm}\par \noindent
\marginpar{\textsc{ToDo}}
\framebox{\begin{minipage}[c]{0.95 \textwidth}
\tt #1 \end{minipage}}\vspace{5 mm}\par}
\makeatletter
\@addtoreset{figure}{section}
\def\thefigure{\thesection.\@arabic\c@figure}
\def\fps@figure{h, t}
\@addtoreset{table}{bsection}
\def\thetable{\thesection.\@arabic\c@table}
\def\fps@table{h, t}
\@addtoreset{equation}{section}
\def\theequation{\thesection.\arabic{equation}}
\makeatother

\allowdisplaybreaks
\def\intprod{\mathbin{\hbox to 6pt{%
                    \vrule height0.4pt width5pt depth0pt
                    \kern-.4pt
                    \vrule height6pt width0.4pt depth0pt\hss}}}

\title{Variational Methods, Multisymplectic Geometry and Continuum
Mechanics}

\author{Jerrold E. Marsden\thanks{
Research partially supported by the
California Institute of Technology
and NSF contract  KDI/ATM-9873133.}
\\{\small Control and Dynamical Systems 107-81}
\\ {\small California Institute of Technology}
\\ {\small Pasadena, CA 91125}
\\ {\small marsden@cds.caltech.edu}
        \and
Sergey Pekarsky$^*$
\\{\small Control and Dynamical Systems 107-81}
\\ {\small California Institute of Technology}
\\ {\small Pasadena, CA 91125}
\\ {\small sergey@cds.caltech.edu}
\and
Steve Shkoller\thanks{
Research partially supported by
NSF-KDI/ATM-9873133 and the Alfred P. Sloan Research Fellowship.}
\\{\small Department of Mathematics}
\\ {\small University of California, Davis}
\\{\small Davis, CA 95616-8633}
\\ {\small shkoller@math.ucdavis.edu}
\and
Matthew West$^*$
\\{\small Control and Dynamical Systems 107-81}
\\ {\small California Institute of Technology}
\\ {\small Pasadena, CA 91125}
\\ {\small mwest@cds.caltech.edu}
}
       \date{\small August 31, 1999; this version April 13, 2000\\}
       \maketitle



\begin{abstract}
This paper presents a variational and multisymplectic formulation of
both compressible and incompressible models of continuum
mechanics on general Riemannian manifolds.
A general formalism is developed
for  non-relativistic first-order
multisymplectic field theories with constraints, such as the
incompressibility constraint. The results obtained in this paper set the
stage for multisymplectic reduction and for the further development
of  Veselov-type multisymplectic discretizations and numerical
algorithms. The latter will be the subject of a companion paper.
\end{abstract}

\maketitle

\tableofcontents

\section{Introduction}

The purpose of this paper is to give a variational multisymplectic
formulation of continuum mechanics from a point of view that will
facilitate the development of a corresponding discrete theory, as
in the PDE Veselov formulation due to Marsden, Patrick, and Shkoller
[1998]. This discrete theory and its relation to finite element
methods will be the subject of a companion paper, Marsden,
Pekarsky, Shkoller, and West [2000].

In this paper, we restrict our attention to
non-relativistic theories, but  on general Riemannian manifolds. The
relativistic case was considered in Kijowski and Tulczyjew [1979],
where the authors take an alternative approach of inverse fields,
effectively exchanging the base and fiber spaces. See also
Fern\'{a}ndez, Garc\'{i}a and Rodrigo [1999].\footnote{There are a
number of reasons, both functional analytic and geometric for
motivating a formulation in terms of {\it direct particle placement
fields} rather than on {\it inverse fields}. For example, in the
infinite dimensional context, this is the setting in which one has
the deeper geometric and analytical properties of the Euler
equations and related field theories, as in Arnold [1966], Ebin and
Marsden [1970], Shkoller [1998], and Marsden, Ratiu and Shkoller [2000].
Moreover, the relativistic approach adopted in Kijowski and Tulczyjew [1979]
  cannot describe an incompressible
fluid or elasticity because the notion of incompressibility
is not defined in the relativistic context.}


Two main applications of our theory are considered---fluid dynamics
and elasticity---each specified by a particular choice of the
Lagrangian density. The resulting Euler-Lagrange equations can be
written in a well-known form by introducing the pressure function
$P$ and the Piola-Kirchhoff stress tensor ${\mathcal P}$ 
(equations (\ref{EL_bf}) and (\ref{EL_el}) below, respectively).

We only consider \emph{ideal}, that is nonviscous, fluid dynamics
in this paper, both compressible and incompressible cases.
In the compressible case, we work out the details for \emph{barotropic}
fluids for which the stored energy is a function of the density.
These results can be trivially extended to \emph{isentropic}
(compressible) fluids, when the stored energy depends also
on the entropy. Both the density and the entropy are assumed
to be some given functions in material representation, so that our
formalism naturally includes \emph{inhomogeneous} ideal fluids. However,
in our discussion of symmetries and corresponding
conservation laws considered in \S\ref{Symm_sec} we restrict ourselves,
for simplicity only, to fluids that are homogeneous in the reference
configuration. We elaborate on this point below.

For the theory of elasticity we restrict our attention to
\emph{hyperelastic} materials, that is to materials whose constitutive
law is derived from a stored energy function. Similarly, we assume
that the material density is some given function which describes a
possibly \emph{heterogeneous} hyperelastic material.

A general formalism for treating constrained multisymplectic
theories is developed in \S\ref{Const_sec}. Often, constraints
that are treated in the multisymplectic context are dynamically
invariant, as with the constraint $\operatorname{div} \mathbf{E} =
0$ in electromagnetism (see, for example, Gotay, Isenberg,
and Marsden [1997]), or $\operatorname{div} \mathbf{E} = \rho$ for
electromagnetism interacting with charged matter. Our main example
of a constraint in this paper is the incompressibility constraint
in fluids, which, when viewed in the standard {\it Eulerian, or
spatial} view of fluid mechanics is often considered to be a {\it
nonlocal} constraint (because the pressure is determined by an
elliptic equation and, correspondingly, the sound speed is
infinite), so it is interesting how it is handled in the
multisymplectic context, which is, by nature, a local formalism.

In the current work we restrict our attention to first-order
theories, in which both the Lagrangian and the constraints depend
only on {\it first} derivatives of the fields. Moreover, we assume
that time derivatives do not enter the constraints, which
corresponds, using a chosen space-time splitting, to holonomic
constraints on the corresponding infinite-dimensional configuration
manifold in {\it material} representation. We briefly discuss the
issues related to  extending this approach to non-holonomic
constraints and to space-time covariant field theories in the last
section.

Symmetries and corresponding momentum maps and conservation theorems
are considered in a separate section (\S\ref{Symm_sec}) since they are
very different for different models of a continuous media, e.g.
homogeneous fluid dynamics has a huge  symmetry, namely the particle relabeling
symmetry, while standard elasticity (usually assumed to be inhomogeneous)
has much smaller symmetry groups, such
as rotations and translations in the Euclidean case.
We emphasize that although the rest of the paper describes general
heterogeneous continuous media, the results of \S\ref{Symm_subsec} only
apply to fluid dynamics that is homogeneous in the reference configuration
(e.g., the fluid starts out, but need not remain, homogeneous), where the
symmetry group is the full group of volume-preserving diffeomorphisms
${\mathcal D}_\mu$. However, these results can be generalized to
inhomogeneous fluids,  in which case the symmetry group is a
\emph{subgroup}
${\mathcal D}^{\rho}_\mu \subset {\mathcal D}_\mu$ that preserves the level
sets of the material density for barotropic fluids, or a \emph{subgroup}
${\mathcal D}^{\rho, \text{ent}}_\mu \subset {\mathcal D}_\mu$ that preserves
the level sets of the material density and entropy for isentropic fluids.
This puts us in the realm of a multisymplectic version of the
Euler-Poincar\'{e} theory -- one needs to introduce additional
\emph{advected quantities} as basic fields to handle this situation
(see the discussion on symmetry and reduction in \S\ref{Disc_sec}).
We remark also that  all continuum mechanics models should satisfy material
frame indifference principle, which, as is well known, can be readily
accomplished by requiring the stored energy function to be a function
of the Cauchy-Green tensor alone (see, e.g. Marsden and Hughes [1983],
Lu and Papadopoulos [1999]).

We finally remark on the notation. The reader is probably aware
that typical fluids and elasticity literatures adhere to completely
different sets of notations, which both differ substantially from those
adopted in multisymplectic theories.
In our notations we follow Gotay, Isenberg, and Marsden [1997].
The companion paper Marsden, Pekarsky, Shkoller, and West [2000]
uses primarily notation from Marsden and  Hughes [1983]
and concentrates on models of continuum mechanics in Euclidean
spaces and their variational discretizations.

\section{Compressible Continuum Mechanics}
\label{Main_sec}

To describe the multisymplectic framework of
continuum mechanics, we must first specify the
covariant configuration  and phase spaces. Once we obtain a
better understanding  of the geometry of these manifolds
we can consider
the dynamics determined by a particular covariant Lagrangian.

\subsection{Configuration and Phase Spaces}

\paragraph{The Jet Bundle.}
We shall set up a formalism in
which a continuous medium is described using sections of a fiber
bundle $Y$ over $X$; here $X$ is the base manifold and
$Y$ consists of fibers $Y_x$ at each point $x \in X$. Sections
of the bundle $\pi_{XY}: Y \rightarrow X $ represent
{\bfi configurations}, e.g. particle placement fields or deformations.

Let $(B, G)$ be a smooth $n$-dimensional compact oriented Riemannian
manifold with a smooth boundary and let $(M, g)$ be a smooth $N$-dimensional
compact oriented Riemannian manifold.
For the non-relativistic case, the base manifold can be chosen to
be a spacetime manifold represented by the product
$X = B \times \mathbb{R} $ of the manifold $B$ together with
time; $(x, t) \in X$. Let us set $x^0 = t$, so that
$x^\mu = (x^i, x^0) = (x^i, t)$, with $\mu = 0, \dots, n, \ i =
1, \dots, n$, denote coordinates on  the $(n+1)$-dimensional
manifold $X$. Construct a trivial bundle $Y$ over $X$ with $M$
being a fiber at each point; that is, $Y = X \times M \ni
(x,t,y)$ with
$y \in M$ --- the fiber coordinate.  This bundle,
$$
\pi_{XY} : Y \rightarrow X; \quad (x,t,y) \mapsto (x,t),
$$
with $\pi_{XY}$ being the projection on the first factor, is the
covariant configuration manifold for our theory.
Let $\mathcal{C}  \equiv C^{\infty}(Y)$ be the
set of smooth sections of $Y$. Then, a section $\phi$ of $
\mathcal{C} $ represents a time dependent configuration.

Let $y^a, \ i = 1, \dots, N$ denote fiber coordinates so
that a section $\phi$ has a coordinate representation
$\phi(x) = (x^\mu, \phi^a(x)) = (x^\mu, y^a)$.
The first jet bundle $J^1 Y$ is the affine bundle over $Y$ whose
fiber above $y \in Y_x$ consists of those linear
maps $\gamma : T_xX \to T_yY$ satisfying $T\pi_{XY} \circ \gamma =
Id_{T_xX}$. Coordinates on
$J^1 Y$ are denoted $\gamma = (x^\mu, y^a, {v^a}_\mu)$.
For a section $\phi$, its tangent map at $x \in X$, denoted
$T_x \phi$, is an element of $J^1 Y_{\phi(x)}$.
Thus, the map $x \mapsto T_x \phi$ is a local section of $J^1 Y$
regarded as a bundle over $X$. This section is denoted $j^1 \phi$
and is called the first jet extension of $\phi$.
In coordinates, $j^1 \phi$ is given by
$(x^\mu, \phi^a(x), \partial_\mu \phi^a)$, where
$\partial_0 \phi^a = \partial_t \phi^a$ and
$\partial_k \phi^a =  \partial_{x^k} \phi^a$.

Notice that we have introduced
\emph{two different} Riemannian structures on the  configuration
bundle.  The internal metric on the spatial
part $B$  of the base manifold $X$ is denoted by $G$ and the fiber, or
field, metric on $M$ is denoted by $g$.
There are two main cases, which we consider in this paper:
\begin{itemize}
\item[(i)] fluid dynamics on a fixed background with fixed boundaries,
when $B$ and $M$ are the same and the fiber
metric $g$ coincides with the base metric $G$;
a special case of this is fluid dynamics on a region in Euclidean
space;
\item[(ii)] elasticity  on a fixed background, when the metric
spaces $(B, G)$ and $(M, g)$ are essentially different.
\end{itemize}
Both approaches result in \emph{background theories}.
The case of relativistic fluid and elasticity was considered
by Kijowski (see, e.g. Kijowski and Tulczyjew [1979]).

Define the following function on the first jet bundle:
\begin{equation}
\label{jacobian_def}
J(x,t,y,v) = \det [v] \sqrt{\frac{\det [g(y)]}{\det [G(x)]}}
: J^1 Y \rightarrow \mathbb{R} .
\end{equation}
We shall see later that its pull-back by a section $\phi$
has the interpretation of the
Jacobian of the linear transformation $D \phi_t$.

A very important remark here is that even though in fluid dynamics
metrics $g$ and $G$ coincide, i.e. on each fiber $Y_x$, $g$ is
a copy of $G$, there is no cancellation because the metric
tensors are evaluated at different points. For instance, in
(\ref{jacobian_def}) $g(y)$ does not coincide with $G(x)$
unless $y=x$ or both metrics are constant.
Hence, only for fluid dynamics in Euclidean spaces, can one
trivially raise and lower indices and drop all metric determinants
and derivatives in the expressions in the next sections.

\paragraph{The Dual Jet Bundle.}

Recall that the dual jet bundle $J^1 Y^*$ is an affine bundle
over $Y$ whose fiber at $y \in Y_x$ is the set of affine maps
from $J^1 Y$ to $\Lambda^{n+1} X_x$, where $\Lambda^{n+1} X$
denotes the bundle of $(n+1)$-forms on $X$.
A smooth section of $J^1 Y^*$ is an affine bundle  map of $J^1 Y$
to $\Lambda^{n+1} X$ covering $\pi_{XY}$. Fiber coordinates on
$J^1 Y^*$ are $(\Pi, {p_a}^\mu)$, which correspond to the affine
map given in coordinates by
${v^a}_\mu \mapsto (\Pi+{p_a}^\mu {v^a}_\mu) d^{n+1} x$.

To define canonical forms on $J^1 Y^*$, another description of
the dual bundle is convenient. Let $\Lambda = \Lambda^{n+1} Y$
denote the bundle of $(n+1)$-forms on $Y$, with fiber over $y \in Y$
denoted by $\Lambda_y$ and with projection
$\pi_{Y \Lambda} : \Lambda \rightarrow Y$.
Let $Z$ be its ``vertically invariant'' subbundle whose fiber is given
by
$$
Z_y = \{ z \in \Lambda_y \mid v \intprod w \intprod z = 0
\ \text{for all} \ v, w \in V_y Y \},
$$
where $V_y Y = \{ v \in T_y Y \mid T \pi_{XY} \cdot v = 0 \}$
is a vertical subbundle.
   Elements of $Z$ can be written uniquely as
$$
z = \Pi d^{n+1} x + {p_a}^\mu d y^a \wedge d^n x_\mu
$$
where $d^n x_\mu = \partial_\mu d^{n+1} x$, so that
$(x^\mu, y^a,  \Pi, {p_a}^\mu)$ give coordinates
on $Z$.

Equating the coordinates $(x^\mu, y^a, \Pi, {p_a}^\mu)$ of $Z$ and of
$J^1 Y^*$ defines a vector bundle isomorphism $Z \leftrightarrow J^1 Y^*$.
This isomorphism can also be defined intrinsically (see
Gotay, Isenberg, and Marsden [1997]).

Define the canonical $(n+1)$-form $\Theta_\Lambda$
on $\Lambda$ by $\Theta_\Lambda (z) = (\pi_{Y \Lambda}^* z)$,
where $z \in \Lambda$. 
The canonical $(n+2)$-form is given by
$\Omega_\Lambda = - d \Theta_\Lambda$.
If $i_{\Lambda Z} : Z \rightarrow \Lambda$ denotes the inclusion,
the corresponding canonical forms on $Z$ are given by
$\Theta = i_{\Lambda Z}^* \Theta_\Lambda$ and
$\Omega  = - d \Theta = i_{\Lambda Z}^* \Omega_\Lambda$.
In  coordinates they have the following representation
$$
\Theta =  {p_a}^\mu d y^a \wedge d^n x_\mu
+ \Pi d^{n+1} x, \quad
\Omega =  d y^a \wedge d {p_a}^\mu \wedge d^n x_\mu
- d \Pi \wedge d^{n+1} x.
$$

\subsubsection*{Ideal Fluid}

We now recall the classical material and spatial descriptions of
ideal (i.e., nonviscous) fluids moving in a fixed region,
i.e., with fixed boundary conditions.
We set $B = M$ and call it the {\bfi reference fluid container}.
A fluid flow is given by a family of 
diffeomorphisms $\eta_t : M \rightarrow M$ with
$\eta_0 = \operatorname{Id}$, where $\eta_t(M)$
is the fluid configuration at some later time $t$. Let $x \in
M$ denote the original position of a fluid particle, then $y
\equiv \eta_t (x) \in M$ is its position at time $t$;
$x$ and $y$ are called {\bfi material} and {\bfi spatial}
points, respectively.
The {\bfi material velocity} is defined by
$V(x,t) = (\partial /\partial t) \eta_t (x)$.
The same velocity viewed as a function of $(y,t)$ is called the
{\bfi spatial velocity}. It is denoted by $u$; that is, $u(y,t)
= V (x(y),t)$, where $x = \eta^{-1}_t (y)$, so that
$u = V \circ \eta^{-1}_t = \dot{\eta} \circ \eta^{-1}_t$.

Thus, in the bundle picture above, the spatial part of the base
manifold $B \subset X$ has the interpretation of the
reference configuration, while an extra dimension of $X$
corresponds to the time evolution. All later configurations
of the fluid are captured by a section $\phi$ of the bundle $Y$,
which gets the interpretation of a particle placement field.
Pointwise this implies that $x$ in the base point $(x,t)$
represents the material point, while $y \in Y_{(x,t)}$ represents
the spatial point and corresponds to a position
$y = \phi(x,t) = \eta_t (x)$ of the fluid particle $x$ at time $t$.

\subsubsection*{Elasticity}

For the theory of elasticity (as well as for fluids with a free
boundary), the base and fiber spaces are generally different; $(B,
G)$ is traditionally called the {\bfi reference configuration},
while
$(M, g)$ denotes the ambient space. For classical $2$ or $3$
dimensional elasticity, $M$ and $B$ have the same dimension,
while for rods and shells models the dimension of the
reference configuration $B$ is less than that of the
ambient space.

For a fixed time $t$, sections of the bundle $Y$, denoted by $\phi_t$,
play the role of {\bfi deformations}; they map reference
configuration $B$ into spatial configuration $M$.
Upon restriction to the space of first jets,
the fiber coordinates $v$ of $\gamma = (x,y, v) \in J^1 Y$ become
partial derivatives $\partial \phi^a / \partial x^\mu$;
they consist of the time derivative of the deformation
$\dot{\phi^a}$ and the {\bfi deformation gradient},
${F^a}_{ i} = \partial \phi^a / \partial x^i$.
The first jet of a section $\phi$ then has the following local
representation $j^1 \phi = ((x,t),\phi(x,t), \dot{\phi}(x,t), F(x,t))
: X \rightarrow J^1 Y$.

Using the map $\phi$, one can pull back and push forward
metrics on the base and fiber manifolds. In particular,
a pull-back of the field metric $g$ on $M$ to $B \subset X$
defines the {\bfi Green deformation tensor} (also called the right
Cauchy-Green tensor) $C$ by $C^\flat = \phi_t^* (g)$,
while a push-forward of the base metric $G$ on $B \subset X$ to $M$
defines the inverse of the {\bfi Finger deformation tensor} $b$
(also called the left Cauchy-Green tensor):
$c = b^{-1} = (\phi_t)_* (G)$.
In coordinates,
\begin{equation}
\label{G&F_def}
C_{ij} (x,t) = g_{ab} {F^a}_{i} {F^b}_{j} (x,t),
\qquad
c_{ab} (y) = G_{ij} {(F^{-1})^i}_{a} {(F^{-1})^j}_{b} (y),
\end{equation}
where $F^{-1}$ is thought of as a function of $y$.
We remark that $C$ is defined whether or not the deformation
is regular, while $b$ and $c$ rely on the regularity of $\phi_t$.
Another important remark is that operations flat $\flat$ and sharp
$\sharp$ are taken with respect to the corresponding metrics on
the space, so that, e.g $(\phi_t^* g)^\sharp \ne \phi_t^* (g^\sharp)$.

Notice that $J$ restricted to the first jets of sections is
the {\bfi Jacobian} of $D \phi_t$, that is, the determinant of the linear
transformation $D \phi_t$; it is given in coordinates by
$$
J (j^1 \phi) = \det [F]
\sqrt{\frac{\det [g]}{\det [G]}} (j^1 \phi)
\ : \ X \rightarrow \mathbb{R} .
$$
It is a scalar function of $x$ and $t$, invariant under coordinate
transformations. Notice also that $J(x,t) > 0$ for regular deformations
with $\phi (x,0) = x, F (x,0) = \operatorname{Id}$ because $J(x,0) = 1$.

\subsection{Lagrangian Dynamics.}
\label{LDyn_subsec}

To obtain the Euler-Lagrange equations
for a particular model of a continuous medium, one needs
to specify a Lagrangian density $\mathcal{L} $. Naturally, it should
contain terms corresponding to the kinetic energy  and to the potential
energy of the medium. Such terms depend on material
properties such as mass density $\rho$ 
as well as on the constitutive relation. The latter is determined
by the form of the potential energy of the material.
We remark that such an approach excludes from our consideration
non-hyperelastic materials whose constitutive laws
cannot be obtained from a potential energy function.

\paragraph{Lagrangian Density.}
Let the mass density $\rho : B \rightarrow \mathbb{R} $ be
given for a particular model of continuum mechanics.
The Lagrangian density
$\mathcal{L}  : \ J^1 Y \rightarrow \Lambda^{n+1} X$
for a multisymplectic model of continuum
mechanics is a smooth bundle map
\begin{multline}
\label{Lag_def}
\mathcal{L}  (\gamma) = L (\gamma)  d^{n+1} x
= {\mathbb K} - {\mathbb P}
= \frac{1}{2} \sqrt{\det [G]} \rho (x) g_{a b} {v^a}_{0} {v^b}_{0} d^{n+1}x \\
- \sqrt{\det [G]} \rho (x) W(x, G(x), g(y), {v^a}_{j}) d^{n+1} x,
\end{multline}
where $\gamma \in J^1 Y$ and $W$ is the {\bfi stored energy function}.
The first term in (\ref{Lag_def}) corresponds
to the kinetic energy of the matter when restricted to first jet
extensions, as ${v^a}_0$ becomes the time derivative $\partial_t \phi^a$
of the section $\phi$.
The second term reflects the potential energy and depends on the spatial
derivatives of the fields (upon restriction to first jet extensions),
i.e. on the deformation gradient $F$.

A choice of the stored energy energy function specifies
a particular model of a continuous medium. While different general
functional forms distinguish various broad classes of materials
(elastic, fluid, etc.), the specific functional forms determine specific
materials. Typically, for elasticity, $W$ depends on the field's partial
derivatives through the (Green) deformation tensor $C$, while for
Newtonian fluid  dynamics, $W$ is
only a function of the Jacobian $J$ (\ref{jacobian_def}).

\paragraph{Legendre Transformations.}
The Lagrangian density (\ref{Lag_def}) determines the
Legendre transformation
${\mathbb F} \mathcal{L}  : J^1 Y \rightarrow J^1 Y^*$.
The conjugate momenta are given by the following expressions
\begin{equation}
\label{Legendre_def}
{p_a}^{0}  = \dfrac{\partial L}{\partial {v^a}_{0}}
= \rho g_{ab} {v^b}_{0} \sqrt{\det [G]},
\qquad
{p_a}^{j}  = \dfrac{\partial L}{\partial {v^a}_{j}}
= -\rho \dfrac{\partial W}{\partial {v^a}_{j}}\sqrt{\det [G]},
\end{equation}
$$
\Pi = L - \dfrac{\partial L}{\partial {v^a}_{\mu}} {v^a}_{\mu}
= \left[-\frac{1}{2}  g_{a b} {v^a}_{0} {v^b}_{0}
- W + \frac{\partial W}{\partial {v^a}_{j}}{v^a}_{j}\right]
\rho \sqrt{\det[G]}
$$
Define the {\bfi energy density} $e$ by
\begin{equation}
\label{E_def}
e = {p_a}^0 {v^a}_0 - L
\quad \operatorname{or, equivalently} \quad
e d^{n+1} x = {\mathbb K} + {\mathbb P},
\end{equation}
then
$$
\Pi = -{p_a}^j {v^a}_j - \sqrt{\det [G]} e.
$$

\paragraph{The Cartan Form.}
Using the Legendre transformation (\ref{Legendre_def}), we can pull-back
the canonical $(n+1)$-form from the dual bundle. The resulting form
on $J^1 Y$ is called the Cartan form and is given by
\begin{multline}
\label{Cartan_def}
\Theta_\mathcal{L}  = \rho g_{ab} {v^b}_0 \sqrt{\det [G]}
dy^a \wedge d^n x_0
-\rho \dfrac{\partial W}{\partial {v^a}_j}\sqrt{\det [G]}
dy^a \wedge d^n x_j \\
+ \left[-\frac{1}{2}  g_{a b} {v^a}_0 {v^b}_0
- W +  \frac{\partial W}{\partial {v^a}_j}{v^a}_j\right]
\rho \sqrt{\det[G]} d^{n+1}x.
\end{multline}
We set $\Omega_\mathcal{L}  = - d \Theta_\mathcal{L} $.

Theorem \ref{thm_cartan} below provides a nicer method for obtaining
the Cartan form via the Calculus of Variations and remains
entirely on the Lagrangian bundle $J^1 Y$.
Moreover, the variational approach is essential for the Veselov
type discretization of our multisymplectic theory.
We present it here for the benefit of the reader,
but remark that it is not essential for our current exposition
and can be omitted on a first reading (see Marsden, Patrick, and Shkoller
[1998] for details).

\paragraph{Variational Approach.}

To make the variational derivation of the equations of motion rigorous
as well as that of the geometric objects, such as the multisymplectic
form and the Noether current, we need to introduce some new notations
(see Marsden, Patrick, and Shkoller [1998]).
These are generalizations of the notations used in the rest of the paper.
They only apply to the variational derivation
described here and later in Subsection
\ref{Symm_subsec} and do not influence the formalism and results
in the rest of the paper. The reason for such generalizations is
very important yet subtle: one should allow for \emph{arbitrary}
and not only \emph{vertical} variations of the sections.

Vertical variations are confined to the vertical subbundle
$VY \subset TY$, $V_y Y = \{ {\mathcal V} \in T_y Y |
T \pi_{XY} \cdot {\mathcal V} = 0 \}$; this  allows
only for \emph{fiber-preserving} variations, i.e., if $\phi (X) \in Y_x$ and
$\tilde{\phi}$ is a new section, then $\tilde{\phi} \in Y_x$.
In general, one should allow for arbitrary  variations in $TY$,
when $\tilde{\phi} \in Y_{\tilde{x}}$ for some $\tilde{x} \ne x$.
Introducing a splitting of the tangent bundle into a vertical and
a horizontal parts, $T_y Y = V_y Y \oplus H_y Y$
($H_y Y$ is not uniquely defined), one can
decompose a general variation into a vertical and horizontal
components, respectively.

Explicit calculations show (see Marsden, Pekarsky, Shkoller,
and West [2000]) that while both vertical and arbitrary variations
result in the same  Euler-Lagrange equations, the Cartan form
obtained from the vertical variations only is missing one term
(corresponding to the $d^{n+1} x$ from on $X$); the horizontal
variations account precisely for this extra term and make
the Cartan form complete.

One can account for general variations
either by introducing new ``tilted sections'', or by introducing
some true new sections that compensate for the horizontal variation.
The later can be implemented in the following way.
Let $U \subset X$ be a smooth manifold with smooth closed boundary.
Define the set of smooth maps
$$
{\mathcal C} = \{ \varphi : U \rightarrow Y |
\pi_{XY} \circ \varphi : U \rightarrow X \text{is an embedding} \}.
$$
For each $\varphi \in {\mathcal C}$, set $\varphi_X=\pi_{XY} \circ \varphi$
and $U_X = \pi_{XY} \circ \varphi (U)$, so that
$\varphi_X : U \rightarrow U_X$ is a diffeomorphism and
$\varphi \circ \varphi_X^{-1}$ is a section of $Y$.
The {\bfi tangent space} to the manifold ${\mathcal C}$ at
a point $\varphi$ is the set $T_{\varphi} {\mathcal C}$ defined by
$$
\{ {\mathcal V} \in C^\infty (X, TY) |
\pi_{Y,TY} \circ {\mathcal V} = \varphi \ \& \
T \pi_{XY} \circ {\mathcal V} = {\mathcal V}_X,
\text{a vector field on} X \}.
$$

Arbitrary (i.e., including both vertical and horizontal) variations
of sections of $Y$ can be induced by a family of maps $\varphi$ defined
through the action of some Lie group. Let ${\mathcal G}$ be a Lie
group of $\pi_{XY}$ bundle  automorphisms $\eta_Y$ covering
diffeomorphisms $\eta_X$. Define the {\bfi action} of ${\mathcal G}$
on ${\mathcal C}$ by composition:
$\eta_Y \cdot \varphi = \eta_Y \circ \varphi$. Hence, while
$\varphi \circ \varphi_X^{-1}$ is a section of $\pi_{U_X,Y}$,
$\eta_Y \cdot \varphi$ induces a section
$\eta_Y \circ (\varphi \cdot \varphi_X^{-1} ) \circ \eta^{-1}_X$
of $\pi_{\eta_X(U_X),Y}$.

A one parameter family of variations can be obtain in the following
way. Let $\varepsilon \mapsto \eta^\varepsilon_Y$ be an arbitrary
smooth path in ${\mathcal G}$ with $\eta^0_Y = e$, and let
${\mathcal V} \in T_\varphi {\mathcal C}$ be given by
$$
{\mathcal V} =  \left. \frac{d}{d \varepsilon} \right|_{\varepsilon=0}
\eta^\varepsilon_Y \cdot \varphi.
$$
Define the {\bfi action function}
$$
S (\varphi) = \int_{U_X}
\mathcal{L}  (j^1 (\varphi \circ \varphi_X^{-1}) )
: \mathcal{C}  \rightarrow \mathbb{R},
$$
and call $\varphi$ a {\bfi critical point (extremum)} of $S$ if
$$
d S (\varphi) \cdot {\mathcal V} \equiv
\left. \frac{d}{d \varepsilon} \right|_{\varepsilon=0}
S (\eta^\varepsilon_Y \cdot \varphi) = 0.
$$

The Euler-Lagrange equations and the Cartan form can be obtained
by analyzing this condition. We summarize the results in the
following theorem  from Marsden, Patrick, and Shkoller [1998]
which illustrates the application of the variational principle to
multisymplectic field theory.

\begin{thm}
\label{thm_cartan}
Given a smooth Lagrangian density $\mathcal{L} :J^1 Y \rightarrow
\Lambda^{n+1}(X)$, there exist a unique smooth section $\ \ D_{EL}
\mathcal{L}  \in C^\infty (Y^{''}, \Lambda^{n+1}(X)
\otimes T^*Y))$ ($Y^{''}$ being the space of second jets of sections)
and a unique differential form $\Theta_\mathcal{L}  \in
\Lambda^{n+1}(J^1 Y)$ such that for any
$\mathcal{V}  \in T_\phi\mathcal{C} $,
and any open subset $U_X$ such that $\overline{U}_X \cap \partial X
= \emptyset$,
\begin{equation}
d S (\varphi) \cdot \mathcal{V}  =  \int_{U_X}
D_{EL} \mathcal{L} (j^2(\varphi \circ \varphi_X^{-1})) \cdot \mathcal{V}
+ \int_{\partial U_X} j^1(\varphi \circ \varphi_X^{-1})^*
[{j^1(\mathcal{V} )} \intprod \Theta_\mathcal{L} ].
\label{s21}
\end{equation}
Furthermore,
\begin{equation}
D_{EL} \mathcal{L} (j^2(\varphi \circ \varphi_X^{-1})) \cdot \mathcal{V}
= j^1(\varphi \circ \varphi_X^{-1})^*[j^1(\mathcal{V} )
\intprod \Omega_\mathcal{L} ] \
\text{ in } U_X.
\label{s21b}
\end{equation}
In coordinates, the action of the Euler-Lagrange derivative
$D_{EL}\mathcal{L} $ on $Y^{''}$ is given by
\begin{multline}
D_{EL} \mathcal{L}  (j^2(\varphi \circ \varphi_X^{-1})) =
\left[
\frac{\partial L }{\partial y^a} (j^1 (\varphi \circ \varphi_X^{-1}))
- \frac{\partial^2 L }{\partial x^\mu \partial {v^a}_\mu}
(j^1 (\varphi \circ \varphi_X^{-1})) \right. \\
- \frac{\partial^2 L }{\partial y^b \partial {v^a}_\mu}
(j^1 (\varphi \circ \varphi_X^{-1}))
\cdot  (\varphi \circ \varphi_X^{-1})^b_{,\mu} \\
- \left. \frac{\partial^2 L }{\partial {v^b}_\nu \partial {v^a}_\mu}
   (j^1 (\varphi \circ \varphi_X^{-1}))
\cdot  (\varphi \circ \varphi_X^{-1})^b_{,\mu \nu} \right]
dy^a \wedge d^{n+1}x ,
\label{s21c}
\end{multline}
while the form $\Theta_\mathcal{L} $ matches the definition of the Cartan
form obtained via Legendre transformation and has the coordinate expression
\begin{equation}
\Theta_\mathcal{L}  = \frac{\partial L}{\partial {v^a}_\mu} dy^a \wedge
d^nx_\mu
+\left( L - \frac{\partial L}{\partial {v^a}_\mu} {v^a}_\mu \right) d^{n+1}x.
\label{s21d}
\end{equation}
\end{thm}
\begin{cor}
\label{cor4.2}
The ($n+1$)-form $\Theta_\mathcal{L} $ defined by the variational
principle satisfies the relationship
$$\mathcal{L} ({\mathfrak z}) = {\mathfrak z}^* \Theta_\mathcal{L} $$
for all holonomic sections ${\mathfrak z} \in C^\infty(\pi_{X,J^1 Y})$.
\end{cor}

Another important general theorem, which we quote from
Marsden, Patrick, and Shkoller [1998],
is the so-called {\bfi multisymplectic form formula}
\begin{thm}
If $\phi$ is a solution of the Euler-Lagrange equation (\ref{s21c}),
then
\begin{equation}
\label{MSFF}
\int_{\partial U_X} (j^1 (\varphi \circ \varphi_X^{-1}))^*
\left[ j^1 \mathcal{V}  \intprod
   j^1 \mathcal{W}  \intprod \Omega_{\mathcal{L} } \right]
=0
\end{equation}
for any $\mathcal{V} , \mathcal{W} $ which solve the first variation
equations of the Euler-Lagrange equations, i.e. any tangent vectors
to the space of solutions of (\ref{s21c}).
\end{thm}

This result is the multisymplectic analog of the fact that the
time $t$ map of a mechanical system consists of canonical
transformations. See Marsden, Patrick, and Shkoller [1998] for the
proofs.

Finally we remark that in order to obtain vertical variations
we can require $\varphi_X$ (and, hence, $\varphi_X^{-1}$)
to be the identity map on $X$. Then, $\phi = \varphi \circ \varphi_X^{-1}$
becomes a true section of the bundle $Y$.

\paragraph{Euler-Lagrange Equations.}
Treating $(J^1 Y, \Omega_{{\mathcal L}})$ as a multisymplectic manifold,
the Euler-Lagrange equations can be derived from the following condition
on a section $\phi$ of the bundle $Y$:
$$
(j^1 \phi )^* ({\mathcal W} \intprod \Omega_{{\mathcal L}}) = 0,
$$
for any vector field ${\mathcal W}$ on $J^1 Y$
(see Gotay, Isenberg Marsden [1997] for the proof).
This translates to the following familiar expression in coordinates
\begin{equation}
\label{EL_eq}
\frac{\partial L}{\partial y^a} (j^1 \phi)
- \frac{\partial }{\partial x^\mu} \left(
\frac{\partial L}{\partial {v^a}_\mu} (j^1 \phi) \right) = 0,
\end{equation}
which is equivalent to equation (\ref{s21c}).

Substituting the Lagrangian density (\ref{Lag_def}) into equation (\ref{EL_eq})
we obtain the following Euler-Lagrange equation for a continuous
medium:
\begin{multline}
\label{EL_cm}
\rho g_{ab} \left( \frac{D_g \dot{\phi}}{D t} \right)^b
-\frac{1}{\sqrt{\det [G]}} \frac{\partial}{\partial x^k}
\left( \rho \dfrac{\partial W}{\partial {v^a}_k}(j^1 \phi)\sqrt{\det [G]}
\right) = \\
   - \rho \dfrac{\partial W}{\partial g_{bc}} \dfrac{\partial g_{bc}}
{\partial y^a} (j^1 \phi),
\end{multline}
where
$$
\left( \frac{D_g \dot{\phi}}{D t} \right)^b \equiv
\frac{\partial \dot{\phi}^b}{\partial t} + \gamma_{bc}^a \
\dot{\phi}^b \dot{\phi}^c
$$
is the covariant time derivative, which
corresponds to {\bfi material acceleration}, with
$$
\gamma_{ab}^c = \frac{1}{2} g^{cd} \left(
\frac{\partial g_{ad}}{\partial y^b} +
\frac{\partial g_{bd}}{\partial y^a} -
\frac{\partial g_{ab}}{\partial y^d} \right)
$$
being the Christoffel symbols associated with the `field' metric $g$.
We remark that all terms in this equation are functions of $x$ and $t$
and hence have the interpretation of material quantities.

Equation (\ref{EL_cm}) is a PDE to be solved for a section $\phi(x,t)$
for a given type of potential energy $W$.
As the gravity here is treated parametrically, the term on the
right hand side of (\ref{EL_cm}) can be thought of as a derivative with
respect to a parameter, and we can define a multisymplectic
analogue of the Cauchy stress tensor $\sigma$ as follows
\begin{equation}
\label{CST_def}
\sigma^{ab} = \frac{2 \rho}{J} \frac{\partial W}{\partial g_{ab}} (j^1 \phi)
\ : \ X \rightarrow \mathbb{R} ,
\end{equation}
where $J = \det [F] \sqrt{\det [g]/\det [G]}$ is the Jacobian.
Equation (\ref{CST_def}) is known in the elasticity literature as the
Doyle-Ericksen formula (recall that our $\rho$ corresponds to
$\rho_{Ref}$, so that the Jacobian $J$ in the denominator disappears).

Another important remark is that the balance of moment of momentum
$$
\sigma^T = \sigma
$$
follows from definition (\ref{CST_def}) and the symmetry of
the metric tensor $g$.

Finally, in the case of Euclidean manifolds with constant metrics $g$ and $G$,
equation (\ref{EL_cm}) simplifies to
\begin{equation}
\label{EL_cm_Euc}
\rho \frac{\partial^2 \phi_a}{\partial t^2}
= \frac{\partial}{\partial x^k}
\left( \rho \dfrac{\partial W}{\partial {v^a}_k} (j^1 \phi) \right).
\end{equation}

\subsubsection*{Barotropic Fluid}

For standard models of barotropic fluids, the potential energy
of a fluid depends only on the Jacobian of the fluid's ``deformation'',
so that $W = W(J(g, G, v))$.
For a general inhomogeneous barotropic fluid, the material density
is a given function $\rho (x)$.
In material representation,
this formalism also includes the case of isentropic fluids in
which there is a possible dependence on entropy. Since, in that
case, entropy is advected, this dependency in the material
representation is subsumed by the dependency of the stored energy
function on the deformation gradient. \footnote{In spatial representation,
of course one has to introduce the entropy as an independent
variable, but this naturally happens via reduction. See Holm,
Marsden and Ratiu [1998] for related results from the point of
view of the Euler-Poincar\'e theory with advected quantities.}

The Legendre transformation can be thought of as defining the
pressure function $P$. Notice that
$$
{p_a}^i = -\rho \frac{\partial W}{\partial {v^a}_i} \sqrt{\det [G]}
= -\rho \frac{\partial W}{\partial J} \frac{\partial J}{\partial {v^a}_i}
\sqrt{\det [G]} =
-\rho \frac{\partial W}{\partial J} J {(v^{-1})_a}^i \sqrt{\det [G]}
$$
and define the pressure function to be
\begin{equation}
\label{pres_def}
P(\phi, x) = - \rho(x) \frac{\partial W}{\partial J}(j^1 \phi (x))
\ : \ \mathcal{C}  \times X \rightarrow \mathbb{R} .
\end{equation}
Then for a given section $\phi$, $P(\phi) :  X \rightarrow
\mathbb{R} $ has the interpretation of the {\bfi material pressure}
which is a function of the material density.
In this case, the Cauchy stress tensor defined by (\ref{CST_def}) is
proportional to the metric with the coefficient being minus
the pressure itself:
$$
\sigma^{ab}(x) = \frac{2 \rho}{J} \frac{\partial W}{\partial J}
\frac{\partial J}{\partial g_{ab}} (j^1 \phi) =
- \frac{2 P}{J} J \frac{1}{2} g^{ab} (j^1 \phi) = - P(x) \; g^{ab}(y(x)).
$$
We remark that this can be a defining equation for the pressure
from which (\ref{pres_def}) would follow.
With this notation the left hand side of
the Euler-Lagrange equations (\ref{EL_cm}) becomes
\begin{multline}
\label{EL_lhs}
\rho g_{ab} \left( \frac{D_g \dot{\phi}}{D t} \right)^b-
\frac{1}{\sqrt{\det [G]}}
\frac{\partial}{\partial x^k} \left(- P J
\left(\left(\frac{\partial \phi}{\partial x}\right)^{-1} \right)_a^k
\sqrt{\det [G]}
\right) = \\
\rho g_{ab} \left( \frac{D_g \dot{\phi}}{D t} \right)^b +
\frac{\partial P}{\partial x^k} J
\left(\left(\frac{\partial \phi}{\partial x}\right)^{-1} \right)_a^k
+ \frac{P\det \left(\dfrac{\partial \phi}{\partial x}\right)}{\sqrt{\det [G]}}
\left(\left(\frac{\partial \phi}{\partial x}\right)^{-1}\right)_a^k
\frac{\partial \det [g]}{\partial x^k}  \\
+ (I) + (II) = \rho g_{ab} \left( \frac{D_g \dot{\phi}}{D t} \right)^b +
\frac{\partial P}{\partial x^k} J
\left(\left(\frac{\partial \phi}{\partial x}\right)^{-1}\right)_a^k +
\frac{P}{2} J g^{bc} \dfrac{\partial g_{bc}}
{\partial y^a},
\end{multline}
where terms (I) and (II) arise from differentiating $\det [v]$ and
${(v^{-1})_a}^k$ and cancel each other. The right hand side of
(\ref{EL_cm}) is given by
$$
- \rho \dfrac{\partial W}{\partial g_{bc}} \dfrac{\partial g_{bc}}
{\partial y^a}  = - \rho \dfrac{\partial W}{\partial J}
\dfrac{\partial J}{\partial g_{bc}} \dfrac{\partial g_{bc}}
{\partial y^a}  = \frac{P}{2} J g^{bc} \dfrac{\partial g_{bc}}
{\partial y^a} .
$$
Notice that the last term in (\ref{EL_lhs}) and in the equation
above coincide, so that the Euler-Lagrange equations for the barotropic
fluid have the following form
\begin{equation}
\label{EL_bf}
\rho g_{ab} \left( \frac{D_g \dot{\phi}}{D t} \right)^b = -
\frac{\partial P}{\partial x^k} J
\left( \left( \frac{\partial \phi}{\partial x}
\right) ^{-1}\right)_a^k,
\end{equation}
where the pressure depends on the section $\phi$ and
the density $\rho$ and is defined by (\ref{pres_def}).
Both the metric $g_{ab}$ and the Christoffel symbols
$\gamma_{ab}^c$ in the covariant derivative are evaluated at
$y=\phi(x,t)$.

One can re-write (\ref{EL_bf}) introducing the {\bfi spatial density}
$\rho_{\operatorname{sp}}=\rho/J$ and defining the {\bfi spatial pressure}
$p(y)$ by the relation $P(x)= p(y(x)) = p( \phi_t (x)$. This yields
$$
\frac{D_g V}{D t}(x,t) = - \frac{1}{\rho_{\operatorname{sp}}}
\operatorname{grad} p \circ \phi \ (x,t),
$$
where $V=\dot{\phi}$.
Compare this to the equations for incompressible ideal hydrodynamics in
\S\ref{Incm_sec}.

\subsubsection*{Elasticity}

The Legendre transformation defines the first Piola-Kirchhoff stress tensor
${{\mathcal P}_a}^i$. It is given, up to the multiple of $-1/\sqrt{\det [G]}$,
by the matrix of the partial derivatives of the Lagrangian with respect to
the deformation gradient:
\begin{equation}
\label{1PK_def}
{{\mathcal P}_a}^i (\phi,x) =
\rho (x) \frac{\partial W}{\partial {v^a}_{i}} (j^1 \phi (x)),
\end{equation}
and for a given section $\phi$, ${{\mathcal P}_a}^i$ is a  stress
tensor defined on $X$.

Notice that the first Piola-Kirchhoff stress tensor is proportional
to the spatial momenta, ${{\mathcal P}_a}^i = - {p_a}^i/ \sqrt{\det [G]}$.
The coefficient $\sqrt{\det [G]}$ arises from the difference in
the volume forms used in standard and multisymplectic elasticity.
In the former, the Lagrangian density is integrated over a space
area using the volume form $\mu_G = \sqrt{\det [G]} d^n x$
associated with the metric $G$, while in the latter, the integration
is done over the space-time using $d^{n+1} x = d t \wedge d^n x$.
We also remark that though traditionally the first Piola-Kirchhoff
stress tensor is normally taken with both indices up, our choice
is more natural in the sense that it arise from the Lagrange transformation
(\ref{1PK_def}) which relates ${{\mathcal P}_a}^i$ with the spatial momenta.

Using  definitions (\ref{CST_def}) and (\ref{1PK_def}),
we can re-write equation (\ref{EL_cm})
in the following form
\begin{equation}
\label{EL_cm2}
\rho g_{ab} \left( \frac{D_g \dot{\phi}}{D t} \right)^b
= {{{\mathcal P}_a}^i}_{|i} + \gamma_{ac}^b
({{\mathcal P}_b}^j {F^c}_j - J g_{bd} \sigma^{dc} ),
\end{equation}
where we have introduced a {\bfi covariant divergence} according to
$$
{{{\mathcal P}_a}^i}_{|i} = \operatorname{DIV} {\mathcal P} =
   \frac{\partial {{\mathcal P}_a}^i}{\partial x^i}
+ {{\mathcal P}_a}^j \Gamma^k_{jk} - {{\mathcal P}_b}^i \gamma^b_{ac} {F^c}_i.
$$
Here $\Gamma^i_{jk}$ are the Christoffel symbols corresponding to the
base metric $G$ on $B \subset X$ (see, for example, Marsden and Hughes [1993]
for an exposition on covariant derivatives of two-point tensors).

We emphasize that in (\ref{EL_cm2}) there is no a-priori relationship
between the first Piola-Kirchhoff stress tensor and the Cauchy stress
tensor, that is, $W$ has the most general form $W(x, G, g, v)$.
Such a relationship can, however, be derived from the fact that
for standard models of elasticity the stored energy function $W$ depends
on the deformation gradient $F$ (i.e. on $v$) and on the field metric
$g$ only via the Green deformation tensor $C$ given by (\ref{G&F_def}),
that is $W=W(C(v,g))$. Thus, the partial derivatives of $W$ with
respect to $g$ and $v$ are related, and  the following equation
$$
{{\mathcal P}_a}^i = J {(\sigma F^{-1})_a}^i
$$
follows from definitions (\ref{CST_def}) and (\ref{1PK_def}).
This relation immediately follows from the form of the stored energy
function; it recovers the Piola transformation law, which in conventional
elasticity relates the first Piola-Kirchhoff stress tensor and the
Cauchy stress tensor.
Substituting this relation in (\ref{EL_cm2}) one easily notices that
the last term on the right hand side cancels, so that the
Euler-Lagrange equation for the standard elasticity model can be
written in the following covariant form
\begin{equation}
\label{EL_el}
\rho \frac{D_g V}{D t} = \operatorname{DIV} {\mathcal P},
\end{equation}
where $V = \dot{\phi}$.
For elasticity in a Euclidean space, this equation simplifies
and takes a well-known form:
$$
\rho \frac{\partial V^a}{\partial t} =
\frac{\partial {\mathcal P^{ai}}}{\partial x^i}.
$$

\section{Constrained Multisymplectic Field Theories}
\label{Const_sec}

Multisymplectic field theory is a formalism for the construction
of Lagran\-gian field theories. This is to be contrasted with the
formalism in which one takes the view of infinite dimensional
manifolds of fields as configuration spaces. The multisymplectic
view makes explicit use of the fact that many Lagrangian field
theories are local theories, that is, the Lagrangian depends only
pointwise on the values of the fields and their derivatives. In
formulating a constrained multisymplectic theory, we will therefore
only be concerned with the imposition of pointwise constraints
$\Phi (\gamma), \ \gamma \in J^1 Y$,
depending on point values of the fields and their derivatives. In
the current work we also restrict our attention to first-order
theories, in which only first derivatives of the fields are
considered.

Despite the pointwise nature of the Lagrangian
${\mathcal L} (\gamma), \ \gamma \in J^1 Y$, the variational principle
assumes variations of local sections over some region $U \subset X$,
that is, it is the action $S(\phi) = \int_U {\mathcal L} (j^1 \phi)$
as a function of sections that is being minimized.
In order to use the theory
of Lagrange multipliers to impose the constraints, it is therefore
necessary to form a function $\Psi (\phi)$ of local sections
which is defined through point values of the constraint $\Phi (j^1 \phi)$
evaluated at the first jets of sections.
It is then possible, however, to use the pointwise nature of the
Lagrangian and the constraint function to derive a purely local
condition, the Euler-Lagrange equations, for the constrained field
variables. We will make these ideas precise in Section \ref{sec:MFT}.

For holonomic constraints it is well
known that Hamilton's principle constrained to the space of
allowable configurations gives the correct equations of motion.
Hamilton's principle can be naturally extended by either
extremizing over the space of motions satisfying the constraints
(so-called vakonomic mechanics), which is appropriate for optimal
control, but not for dynamics, or by requiring stationarity of the
action with respect to variations which satisfy the constraints (the
Lagrange-d'Alembert or virtual work principle). The equations of
motion derived in each case are, however, different.

Derivations from balance laws (Jalnapurkar [1994]), evidence from
experiments (Lewis and Murray [1995]) and comparison to Gauss' Principle
of Least Constraint and the Gibbs-Appell equations (Lewis [1996])
indicates that it is the Lagrange-d'Alembert principle which gives the
correct equations of motion; see Bloch and Crouch [1999] for
further discussion and references.

While the subject of linear and affine non-holonomic constraints is
relatively well-understood (see Bloch, Krishnaprasad, Marsden, and
Murray [1996]), it is less clear how to proceed for non-linear
non-holonomic constraints. Part of the problem lies in the lack of
examples for which the correct equations are clear from physical
grounds.
In the context of constrained field theories, however, there are many
cases where nonlinear constraints involving spatial derivatives of the
fields need to be applied, such as incompressibility in fluid
mechanics, and it is clear what the physically correct equations
should be. Here we deliberately avoid the use of the term
non-holonomic, to avoid confusion with its standard meaning in the
ODE context, where it applies only to time derivatives. Other
examples of nonlinearly constrained field theories include
constrained director models of elastic rods and shells.

The fact that the constraints involve only spatial and not time
derivatives means that imposing the constraints is equivalent to
restricting the infinite-dimensional configuration manifold used to
formulate the theory as a traditional Hamiltonian or Lagrangian field
theory. In this case, the constraint is simply a holonomic or
configuration constraint and it is known that restricting Hamilton's
principle to the constraint submanifold gives the correct equations
for the system.

\subsection{Lagrange Multipliers}
\label{sec:LM}

The Lagrange multiplier theorem naturally makes use of the dual of the
space of constraints. In a finite-dimensional setting this is a well
defined object, with all definitions being equivalent. When
considering infinite-dimensional constraint spaces, however, the issue
of what is being used as the dual becomes less clear and more
important.

We shall consider constrained multisymplectic field theories for which
the constraint space is the space of smooth sections of a particular
vector bundle. In the case of the incompressibility constraint, the vector
space is one-dimensional and the constraint bundle is, effectively,
the space of real valued functions on the base space $X$.
A dual of the constraint
space is then defined with respect to an inner product structure
on the vector bundle. This is made explicit
in the following statement of the Lagrange multiplier theorem
where we assume that fields and Lagrange multipliers are sufficiently
regular (see Luenberger [1969]).

\begin{thm}[Lagrange multiplier theorem]
Let $\pi_{\mathcal{M} , \mathcal{E} } : \mathcal{E}  \to \mathcal{M} $
be an inner product bundle over a smooth
manifold $\mathcal{M} $, $\Psi$ a smooth section of
$\pi_{\mathcal{M} , \mathcal{E} }$, and $h : \mathcal{M}
\to \mathbb{R}$ a smooth function. Setting $\mathcal{N}  = \Psi^{-1}(0)$,
the following are equivalent:
\begin{enumerate}
\item $\varphi \in \mathcal{N} $ is an extremum of $h|_{\mathcal{N} }$
\item there exists an extremum $\bar{\varphi} \in \mathcal{E} $ of $\bar{h} :
\mathcal{E}  \to \mathbb{R}$ such that $\pi_{\mathcal{M} ,\mathcal{E} }
(\bar{\varphi}) = \varphi$
\end{enumerate}
where $\bar{h}(\bar{\varphi}) = h(\pi_{\mathcal{M} ,\mathcal{E} }
(\bar{\varphi})) -
\langle \bar{\varphi},\Psi(\pi_{\mathcal{M} ,\mathcal{E} }(\bar{\varphi}))
\rangle_{\mathcal{E} }$.
\end{thm}

If $\mathcal{E}$ is a trivial bundle over $\mathcal{M}$, then in
coordinates of the trivialization we have $\bar{\varphi} = (\varphi,
\lambda)$, where  $\lambda : \mathcal{M}  \rightarrow \mathcal{E}
/\mathcal{M}$ is a Lagrange multiplier function.

In the next section we shall use this theorem to relate the
constrained Hamilton's principle with the
extremum of the augmented action integral
which contains the constraint paired with a Lagrange multiplier.
Both of them result in constrained Euler-Lagrange
equations. We shall furthermore demonstrate that, using
the trivialization coordinates,
these equations can be equivalently obtained from
a Lagrangian defined on an extended configuration
bundle. In this picture, the Lagrange multiplier corresponds
to a new field, which extends the dimension of the fiber
space, and the augmented Lagrangian contains an additional
part corresponding to the pairing of this field with the constraint.
The Euler-Lagrange equations of motion then follow from
\emph{unconstrained} Hamilton's principle in a standard way.

\subsection{Multisymplectic Field Theories}
\label{sec:MFT}

In the setting above, the \emph{configuration bundle} is a fiber
bundle $\pi_{X,Y} : Y \to X$ and $\pi_{Y,J^1Y} : J^1Y \to Y$ is the
corresponding first jet bundle with $x^{\mu}$ and $y^a$ being a
local coordinate system on $X$ and $Y$ respectively, and
${v^a}_{\mu}$ the fiber coordinates on $J^1Y$.

Choose the \emph{configuration manifold} $\mathcal{M} $ to be the space
$\mathcal{C} $ of smooth sections $\phi$ of $\pi_{X,Y}$.
Recall that for a Lagrangian density $\mathcal{L}  : J^1Y \to \Lambda^{n+1}X$,
a section $\phi \in \mathcal{M} $ is said to satisfy
Hamilton's principle if $\phi$ is an extremum of
the action function $S (\phi) = \int_X \mathcal{L} (j^1 \phi) :
\mathcal{M}  \to \mathbb{R} $.
Choose the $h$ above to be the action function $S$ and
use $\bar{S}$ instead of $\bar{h}$.


To apply the Lagrange multiplier theorem we need to define
constraints as a section of some bundle
$\mathcal{E}  \rightarrow \mathcal{M} $ (below called the
constraint bundle).
As mentioned above, we restrict our attention to constraints $\Phi$
which depend only on point values of the fields and their derivatives.
Using such constraints we can construct induced constraints $\Psi$
according to (\ref{glob_loc_const}). This is made precise below.
We point out, however, that our treatment excludes inherently global
constraints, such as those on the inverse Laplacian of the field, which
can not be derived from pointwise values.

On the other hand, we also exclude from the consideration a
(simple) case when the constrained subbundle of $J^1Y$ can be
trivially realized as the first jet of some subbundle of $Y$.

Define an inner product vector bundle
$\pi_{X,\mathcal{V} } : \mathcal{V}  \to X$ with the inner product denoted by
$\langle \cdot, \cdot \rangle_{\mathcal{V} }$ whose fibers are
isomorphic to $\mathbb{R} ^n$.
Let $\mathcal{C} ^\infty (\mathcal{V} )$ be the inner product
space of smooth sections of $\pi_{X,\mathcal{V} }$ with the inner
product given by
$$
\langle a,b \rangle = \int_X \langle a(x),b(x) \rangle_{\mathcal{V} }
d^{n+1} x.
$$

The constraint function is an $\mathbb{R} ^n$-valued function on $J^1 Y$:
$$
\Phi : J^1 Y \rightarrow \mathbb{R} ^n.
$$
We say that a point $\gamma \in J^1 Y$ satisfies the constraint
if $\Phi (\gamma) = 0$.
By restricting $\Phi$ to the space of first jets of sections $\phi$ of $Y$,
we can define the {\bfi induced constraint function} $\Psi$
from $\Phi$ by setting
\begin{equation}
\label{glob_loc_const}
\Psi(\phi) (x) = \Phi((j^1 \phi)(x))
\end{equation}
for all $\phi \in \mathcal{M} $ and $x \in X$.
By construction, $\Psi$ is a map from the space $\mathcal{M} $ of
sections of $Y$ to the space $\mathcal{C} ^\infty (\mathcal{V} )$
of sections of $\mathcal{V} $, hence
it can be thought of as a smooth section $\Psi : \mathcal{M}  \rightarrow
\mathcal{E} $ of the {\bfi constraint bundle} $\mathcal{E} $.
This bundle is the trivial inner product bundle given by $\mathcal{M}  \times
\mathcal{C} ^\infty (\mathcal{V} )$ over $\mathcal{M} $.
Then, a configuration $\phi \in \mathcal{M} $ is said to \emph{satisfy
the constraints} if
$\Phi((j^1 \phi)(x)) = 0$ for all $x \in X$, that is, the section
$\Psi(\phi)$ must be a zero function on $X$.
This implies that the space of configurations which
satisfy the constraints is given by $\mathcal{N}  = \Psi^{-1}(0)$.

The \emph{constrained Hamilton's principle} now seeks a
$\phi \in \mathcal{N} $ which is an extremum of $S |_\mathcal{N} $.
The Lagrange multiplier theorem given in the previous subsection
can be applied to conclude that this is
equivalent to the existence of $\bar{\phi} \in \mathcal{E} $ with
$\pi_{\mathcal{M} ,\mathcal{E} } (\bar{\phi}) = \phi$
which is an extremum of $\bar{S}$.
Using the coordinates of the trivialization of $\mathcal{E} $ we can write
$\bar{\phi}=(\phi, \lambda)$, where
$\phi = \pi_{\mathcal{M} ,\mathcal{E} } (\bar{\phi})$ is the base point,
i.e. section $\phi$ of $Y$, and $\lambda$ is thought of as a
section of the bundle $\pi_{X,\mathcal{V} }$, i.e. an
$\mathbb{R} ^n$-valued function on $X$.
Then $\bar{S}  : \mathcal{E}  \to \mathbb{R} $ is given by
\begin{align*}
\bar{S}(\bar{\phi}) &= S(\phi) - \langle \lambda, \Psi(\phi)
\rangle_{\mathcal{E} } \\
&= \int_X {L((j^1 \phi)(x)) d^{n+1} x}
- \int_X{ \langle \lambda(x),\Phi((j^1 \phi)(x)) \rangle_{\mathcal{V}
} d^{n+1} x} \\
&= \int_X\left[
L((j^1 \phi)(x))
- \langle \lambda (x), \Phi((j^1 \phi)(x)) \rangle_{\mathcal{V} }
\right] d^{n+1} x
\end{align*}

In the next section, we demonstrate these constructions for the
incompressibility constraint for continuum theories.

The requirement that $\bar{S}$ be stationary with respect to
variations in $\lambda$ at the point $\bar{\phi}$ implies that
\begin{align*}
0 = \frac{\delta\bar{S}}{\delta\lambda}(\bar{\phi})\cdot
\delta\lambda
&= \int_X\left[
- \langle \delta\lambda (x), \Phi((j^1 \phi)(x)) \rangle_{\mathcal{V} }
\right] d^{n+1} x
\end{align*}
for all variations $\delta \lambda$, and thus that
$\Phi((j^1 \phi)(x)) = 0$ for all $x \in X$. This therefore recovers
the condition that $\phi$ must satisfy the constraints.

Stationarity of $\bar{S}$ with respect to variations in $\phi$
can be used to derive the \emph{constrained Euler-Lagrange equations},
which have the form
\begin{multline}
\label{cons_EL}
\frac{\partial}{\partial x^{\mu}}\left(
\frac{\partial L}{\partial {v^a}_{\mu}}((j^1 \phi)(x))
\right)
- \frac{\partial L}{\partial y^a}((j^1 \phi)(x)) \\
+ \left\langle \lambda (x),
\frac{\partial \Phi}{\partial y^a}((j^1 \phi)(x)) \right\rangle
- \frac{\partial}{\partial x^{\mu}} \left\langle \lambda (x),
\frac{\partial \Phi}{\partial {v^a}_{\mu}}((j^1 \phi)(x))
\right\rangle
= 0.
\end{multline}

Alternatively, one can handle the constraints  by
introducing another bundle, denoted by $E$, which is a product bundle
over $X$ with fibers diffeomorphic to $Y_x \times \mathcal{V} _x$.
One can think of $E$ as a configuration bundle
of the corresponding \emph{unconstrained} system whose
Lagrangian contains an additional part corresponding to the
pairing of the constraint with the Lagrange multiplier:
$$
L_\Phi = L + \langle \lambda , \Phi \rangle_{\mathcal{V} }.
$$
The Euler-Lagrange equations of motion then follow from
\emph{unconstrained} Hamilton's principle in a standard way
and coincide with (\ref{cons_EL}). We work out the details for the
incompressibility constraint in the next section.

\section{Incompressible Continuum Mechanics}
\label{Incm_sec}

In this section we shall consider the incompressibility constraint
using the multisymplectic description of continuum mechanics.
The main issue is a proper interpretation of the constraint using
the Lagrange multiplier formalism developed in the previous section.

\subsection{Configuration and Phase Spaces}

Here we briefly summarize the results.  See the analogous parts of
\S\ref{Main_sec} for more details.

\paragraph{Extended Covariant Configuration Bundle.}
The fibers of $\mathcal{V} $ in this case are one-dimensional and
sections $\bar{\phi} = (\phi, \lambda)$ of $E$  contain
both the deformation field and the Lagrange multiplier, i.e.,
$E$ denotes a bundle over $X$ whose fibers are diffeomorphic to the product
manifold $M \times \mathbb{R} $ with the projection map
$$
\pi_{XE} : E \rightarrow X; \quad (x,t,y,\lambda) \mapsto (x,t).
$$
Here $\lambda$ is a section of the trivial bundle
$X \times \, \mathbb{R} $ over $X$,
which can be thought of as a function $\lambda (x,t)$ on $X$.
The phase space is then the first jet bundle $J^1 E$ with coordinates
$\bar{\gamma} = (x^\mu, y^a, \lambda, {v^a}_\mu, \beta_\mu)$;
the first jet extension of a section $\bar{\phi}=(\phi,\lambda)$
has  the following coordinate representation
$(x^\mu, y^a, \lambda, \partial_\mu \phi^a, \partial_\mu \lambda)$.

\paragraph{The Dual Jet Bundle.}

We can consider the affine dual bundle $J^1 E^*$ as a ``vertically invariant''
subbundle $Z$ of the bundle $\Lambda = \Lambda^{n+1} E$ of all
$(n+1)$-forms on $E$. Elements of $Z$ can be written uniquely as
$$
z = \Pi d^{n+1} x + {p_a}^\mu d y^a \wedge d^n x_\mu
+ \pi^\mu d \lambda \wedge d^n x_\mu,
$$
where $d^n x_\mu = \partial_\mu d^{n+1} x$, so that
$(x^\mu, y^a, \lambda, \Pi, {p_a}^\mu, \pi^\mu)$ give coordinates
on $Z$.

The canonical $(n+1)$-form is constructed in a standard manner
and in the above coordinates has the following representation
$$
\Theta =  {p_a}^\mu d y^a \wedge d^n x_\mu
+ \pi^\mu d \lambda \wedge d^n x_\mu + \Pi d^{n+1} x.
$$
We set $\Omega = - d \Theta$.

The \emph{primary constraint manifold} ${\mathfrak C}$ is a
subbundle of the dual jet bundle and corresponds to the
incompressibility constraint. The pull-back of the inclusion
map $i_{{\mathfrak C}} : {\mathfrak C}
\rightarrow J^1 E^*$ defines a degenerate $(n+2)$-form
$\Omega_{{\mathfrak C}}$ on ${\mathfrak C}$.
We shall discuss the explicit form of the
constraint in the next subsection.

\paragraph{Incompressibility Constraint.}

Recall that such a constraint in, for example, incompressible fluid dynamics,
is a reflection of the divergence-free property of the Eulerian
fluid velocity and, hence, has a pointwise character.
The divergence-free character of the velocity field arises
from the requirement that the particle placement map be volume
preserving at each instant of time.
Then, according to the general
theory of constrained multisymplectic fields outlined above,
it can be obtained from a pointwise constraint $\Phi$ defined
on the first jet bundle $J^1 Y$.

For $\gamma \equiv (x^\mu, y^a, {v^a}_\mu) \in J^1 Y$
we impose the constraint $\Phi(\gamma) = 0$ on the Jacobian of the
deformation, where
\begin{equation}
\label{const_def}
\Phi: J^1 Y  \rightarrow \mathbb{R} ; \ \
\gamma \mapsto J(\gamma) - 1, \ \
J(\gamma) = \det [v] \sqrt{\frac{\det [g(y)]}{\det [G(x)]}},
\end{equation}
where we have used the definition of $J$ given in (\ref{jacobian_def}).
Restricting $\Phi$ to the first jet of a section $\phi$ results in a
constraint on the matrix of spatial partial derivatives
$\partial_j \phi^a$.

For the Lagrange multiplier itself, we choose  the following Ansatz
\begin{equation}
\label{sp_pres_def}
\lambda (x) = \sqrt{\det [G]} P(x) \ : \ X \rightarrow \mathbb{R} ,
\end{equation}
where $P$ will be shown later to have the interpretation of the
\emph{material} pressure. Equation (\ref{sp_pres_def}) guarantees that
$\lambda$ transforms like a density under the transformations
of the base manifold $X$, so that the pairing of $\lambda$
and $\Phi$, defined by integrating over $X$,
has the correct transformation law.

\subsection{Lagrangian Dynamics.}

As we have already mentioned, the main distinguishing feature of
incompressible models of continuum mechanics is the presence of the
constraint  (\ref{const_def}).  We shall now explain how this
modification to the Lagrangian alters the Legendre transform as
well as the Euler-Lagrange equations.

\paragraph{The Lagrangian Density.}
The Lagrangian density $\mathcal{L}  : \ J^1 E \rightarrow
\Lambda^{n+1} X$ for the  multisymplectic model of incompressible
continuum mechanics is a smooth bundle map defined by
\begin{equation}
\label{Lag_incm_def}
\mathcal{L} _\Phi (\bar{\gamma}) = 
\left( L (\gamma) + \lambda \cdot \Phi (\gamma) \right) \ d^{n+1} x =
{\mathbb K} - {\mathbb P} + \lambda \cdot \Phi \; d^{n+1} x,
\end{equation}
where $L$ (i.e. ${\mathbb K}$ and ${\mathbb P}$) is given by
(\ref{Lag_def}) and depends on the choice of the stored energy function $W$.

\paragraph{The Legendre Transformation.}

For the above choice of the Lagrangian, the Legendre transform
thought of as a fiber preserving bundle map
${\mathbb F} \mathcal{L} _\Phi : J^1 E \rightarrow J^1 E^*$
over $E$ is degenerate due to the constrained character of
the dynamics. Indeed,  the Lagrange multiplier $\lambda$ is
a cyclic variable as the Lagrangian (\ref{Lag_incm_def}) does not
depend on its derivatives, $\beta_\mu$. Hence, the conjugate
momentum to
$\lambda$ is identically zero:
$\pi^\mu \equiv \partial L_\Phi / \partial \beta_\mu = 0$.
The set $\{ \pi^\mu = 0 \}$ defines the primary constraint set as
a subset of the dual bundle $J^1 E^*$ to which we restrict the
Legendre transformation to make it non-degenerate.
The rest of the momenta are given by the following expressions
$$
{p_a}^{0}  = \dfrac{\partial L_\Phi}{\partial {v^a}_{0}}
= \rho g_{ab} {v^b}_{0} \sqrt{\det [G]},
$$
\begin{equation}
\label{Legendre_incm_def}
{p_a}^{j}  = \dfrac{\partial L_\Phi}{\partial {v^a}_{j}}
= \left( P J {(v^{-1})_a}^j
- \rho \dfrac{\partial W}{\partial {v^a}_{j}} \right)
   \sqrt{\det [G]},
\end{equation}
$$
\Pi 
= \left[-\rho \frac{1}{2}  g_{a b} {v^a}_{0} {v^b}_{0}
+ \rho \frac{\partial W}{\partial {v^a}_{j}}{v^a}_{j}
- \rho W - P \left(J (n-1) + 1 \right) \right]
\sqrt{\det[G]}.
$$

\paragraph{Euler-Lagrange Equations}

Using the trivialization $(\phi, \lambda)$, we now consider the
Euler-Lagrange equations for a section  $\bar{\phi}$ of $E$,  both
with respect to the deformation $\phi$ and with respect to the
Lagrange multiplier $\lambda$. The former can be written
in coordinates  as follows
\begin{equation}
\label{EL_equiv}
\frac{\partial L_\Phi}{\partial y^a} (j^1 \bar{\phi})
- \frac{\partial }{\partial x^\mu} \left(
\frac{\partial L_\Phi}{\partial {v^a}_\mu} (j^1 \bar{\phi}) \right) = 0.
\end{equation}
The Euler-Lagrange equation for $\lambda$ trivially recovers the
constraint $\Phi = 0$ itself restricted to the first jet :
\begin{equation}
\label{EL_const}
\frac{\partial L_\Phi}{\partial \lambda} (j^1 \bar{\phi})
- \frac{\partial }{\partial x^\mu} \left(
\frac{\partial L_\Phi}{\partial \beta_\mu} (j^1 \bar{\phi}) \right)
= \Phi (j^1 \phi) d^{n+1} x = (J(j^1 \phi) - 1) d^{n+1} x = 0.
\end{equation}
These two equations are to be solved for the Lagrange multiplier
$\lambda$ (equivalently, for the pressure $P$) and for the section
$\phi$.

Substituting Lagrangian (\ref{Lag_incm_def}) into (\ref{EL_equiv}),
we obtain the Euler-Lagrange equation (\ref{EL_cm}) modified by
the pressure term:
\begin{multline}
\label{EL_cm_incm}
\rho g_{ab} \left( \frac{D_g \dot{\phi}}{D t} \right)^b
-\frac{1}{\sqrt{\det [G]}} \frac{\partial}{\partial x^k}
\left( \rho \dfrac{\partial W}{\partial {v^a}_k}(j^1 \phi)
\sqrt{\det [G]} \right) = \\
    - \rho \dfrac{\partial W}{\partial g_{bc}} \dfrac{\partial g_{bc}}
{\partial y^a} (j^1 \phi)
- \dfrac{\partial P}{\partial x^k}
{( v^{-1})_a}^k \; J (j^1 \phi).
\end{multline}
Notice that in the case of parameterized \emph{non-constant} metrics,
the extra pressure term in (\ref{Legendre_incm_def}) gives rise to the term
$$
\frac{\partial}{\partial y^b}
\left( \left( P \; J {(v^{-1})_a}^j  \sqrt{\det [G]} \right) (j^1 \phi)\right)
\ \frac{\partial y^b}{\partial x^j},
$$
which follows from the chain rule applied to $\partial_{x^j} g(y(x))$.
This term exactly cancels another term coming from differentiating
the constraint with respect to $y$ due again to the composition $g=g(y)$ :
$$
\lambda \frac{\partial \Phi}{\partial y^a} =
\frac{P}{2} J g^{bc} \dfrac{\partial g_{bc}}{\partial y^a}\sqrt{\det[G]}
$$
and other cancellations occur as in equation (\ref{EL_lhs}).

In the case of Euclidean manifolds with constant metrics $g$ and $G$,
the Euler-Lagrange equations  simplify to
\begin{equation}
\label{EL_cm_incm_Euc}
\rho \frac{\partial^2 \phi_a}{\partial t^2}
= \frac{\partial}{\partial x^k}
\left( \rho \dfrac{\partial W}{\partial {v^a}_k} (j^1 \phi) \right)
- \dfrac{\partial P}{\partial x^k}
{( v^{-1})_a}^k \; J (j^1 \phi)
\end{equation}
together with the constraint (\ref{EL_const}).

\subsection{Incompressible Ideal Hydrodynamics}

For fluid dynamics, the stored energy term in the Lagrangian
is a constant function precisely because of the incompressibility
constraint. Indeed, as we have mentioned above, $W$ in ideal fluid
models is a function of the Jacobian $J$, but the latter is constrained
to be $1$. For simplicity, consider an  ideal homogeneous incompressible
fluid, so that the material density $\rho$ is constant,
and we can set $\rho=1$ (for inhomogeneous fluids the dependence
of material density on the point $x$ is implicit in the pressure function
$P$).

The Lagrangian is given by (\ref{Lag_incm_def}) with ${\mathbb P} =
\operatorname{const}$. Hence, two terms in  equation
(\ref{EL_cm_incm}) which correspond to the derivatives of $W$
vanish, so that the dynamics of an incompressible ideal fluid is
described by
\begin{equation}
\label{EL_if_incm_a}
g_{ab} \left( \frac{D_g \dot{\phi}}{D t} \right)^b =
- \dfrac{\partial P}{\partial x^k} J
\left(\left(\frac{\partial \phi}{\partial x}\right)^{-1}\right)_a^k,
\end{equation}
together with the constraint
\begin{equation}
\label{EL_if_incm_b}
J (j^1 \phi) =
\left( \frac{\sqrt{\det [G \circ \phi]}}{\sqrt{\det [G]}}
\det \left(\frac{\partial \phi}{\partial x}\right) \right) (x,t) = 1,
\end{equation}
where we have used the fact that $g=G$.

Compare (\ref{EL_if_incm_a}) with (\ref{EL_bf}) and notice that the
incompressibility constraint $J (j^1 \phi) = 1$ implies that the
{\bfi spatial density} $\rho_{\operatorname{sp}} = \rho / J$
is constant, e.g., $1$. Introducing
the {\bfi spatial pressure} $p = P \circ \phi_t ^{-1}$,
the above equation can be written as
\begin{equation}
\label{EL_if_incm2}
\frac{D_g \dot{\phi}}{D t} (x,t) = - \operatorname{grad} p \circ \phi (x,t),
\end{equation}
where we have set $\rho_{\operatorname{sp}} = 1$.
We remark again that the covariant derivative is evaluated
at $y=\phi(x,t)$.

\paragraph{A New Look at the Pressure.}

Here we shall demonstrate that the same equations of motion are
obtained if the potential energy in the Lagrangian (\ref{Lag_incm_def})
is not set to a constant, but rather is treated as a function of the
Jacobian, $W=W(J)$. This will also clarify the relation between the two
definitions of pressure that we have thus far examined.

Recall the definition of the pressure function for barotropic fluids
given by (\ref{pres_def}) as a partial derivative of the stored energy
function $W$ with respect to the Jacobian $J$.
Compare this to the definition (\ref{sp_pres_def}) of the pressure
as a Lagrange multiplier corresponding to the incompressibility
constraint (\ref{const_def}) (modulo a $\sqrt{\det [G]}$ term).
In this subsection we shall denote these objects by $P_W$ and $P_\lambda$,
respectively:
$$
P_W = - \rho \frac{\partial W}{\partial J},
\qquad
P_\lambda = \frac{1}{\sqrt{\det [G]}} \lambda.
$$

The resulting Euler-Lagrange equations can be obtained by combining
(\ref{EL_bf}) with (\ref{EL_if_incm_a}) and are given by:
$$
g_{ab} \left( \frac{D_g \dot{\phi}^b}{D t} \right)^b =
- \dfrac{\partial (P_W + P_\lambda)}{\partial x^k} J
\left(\left(\frac{\partial \phi}{\partial x}\right)^{-1}\right)_a^k
$$
together with the constraint (\ref{EL_if_incm_b}).
We can define a new pressure function
\begin{equation}
\label{new_pres_def}
P = P_W + P_\lambda.
\end{equation}

Notice that when the constraint $J=1$ is enforced by the Euler-Lagrange
equation (\ref{EL_const}), $P_W (J) = \operatorname{const}$,
so that $P = P_{\lambda} + \operatorname{const}$.
This is equivalent to a re-definition of the Lagrange multiplier $\lambda$.
At the same time, the above Euler-Lagrange equation coincides with
(\ref{EL_if_incm_a}) because $\partial_k P = \partial_k P_\lambda$.

\paragraph{Relation to Standard Ideal Hydrodynamics.}

Recall the Lie-Poisson description of fluid dynamics
as a right-invariant system on the group ${\mathcal D}_\mu (M)$
of volume-preserving diffeomorphisms of a Riemannian manifold $(M, G)$.
Here we follow Marsden and Ratiu [1999] and
Arnold and Khesin [1998], using our notations. The Lie algebra
of ${\mathcal D}_\mu (M)$ is the algebra of divergence-free
vector fields on
$M$ tangential to the boundary with minus the Jacobi-Lie
bracket.  The $L^2$ inner-product on this algebra is given by
$$
\langle u, v \rangle _{L^2} = \int_M \langle u(x), v(x)\rangle_G \mu,
$$
where $\mu$ is the Riemannian volume form on $M$.

We extend this metric by right invariance to
the entire group. The resulting Riemannian manifold with right
invariant $L^2$ metric, denoted by
$({\mathcal D}_\mu (M), L^2)$, is the configuration space
for the Lie-Poisson or Euler-Poincar\'e model of ideal
hydrodynamics. Its tangent bundle is the phase space, so that
$(\eta_t, \dot{\eta}_t)$ are the basic ``coordinates''; here
$\eta_t \in {\mathcal D}_\mu (M)$ is a diffeomorphism that
transforms the reference fluid configuration to its
configuration at time $t$. Then, using the kinetic energy of
fluid particles as a Lagrangian, one obtains the following covariant
equations of motion:
\begin{equation}
\label{Euler_on_group_cov}
\frac{D \dot{\eta}}{D t} (x) = - \operatorname{grad} p \circ \eta (x),
\end{equation}
where
$$
\frac{D \dot{\eta}}{D t} = \ddot{\eta} + \Gamma_\eta (\dot{\eta}, \dot{\eta})
$$
denotes covariant material time derivative with respect to the metric
($\Gamma_\eta$ denotes the connection associated to the metric)
and $p$ is the \emph{spatial pressure}.
Notice that covariant derivative is evaluated at $\eta(x)$.

Now define $\eta_t(x) = \eta (t,x)$ to be the flow of the
time-dependent vector field $u(t,x)$, so that
$\partial_t \eta (t,x) = u (t, \eta(t,x))$.
Then composing (\ref{Euler_on_group_cov}) on the right with $\eta^{-1}$
gives the classical Eulerian description of incompressible ideal
fluids :
$$
\partial_t u (t,x) + (u \cdot \nabla ) u = - \operatorname{grad} p,
\qquad \operatorname{div} u = 0.
$$
Taking the divergence of both sides of this expression yields the
equation for the pressure
\begin{equation}
\label{Lapl_pressure}
\Delta p = - \operatorname{div} \left( (u \cdot \nabla ) u \right).
\end{equation}

One readily notices that equations (\ref{EL_if_incm2}) and
(\ref{Euler_on_group_cov}) coincide provided $\eta_t(x) = \phi(x,t)$.
Upon this identification, the Euler-Lagrange equations for
the multisymplectic model of incompressible ideal hydrodynamics recover
the well-known evolution of fluid diffeomorphisms (\ref{Euler_on_group_cov}).
Similarly, taking the divergence of both sides of (\ref{EL_if_incm2})
results in the  Poisson equation on the pressure (\ref{Lapl_pressure}).

\subsection{Incompressible Elasticity}

In a manner similar to the previous subsection,
we modify the elasticity Lagrangian
by the constraint and extend the phase space to include the Lagrange
multiplier. Recall that the stored energy is a function of the
Green deformation tensor $W=W(C)$ and use
the definition of the first Piola-Kirchhoff stress tensor ${{\mathcal P}_a}^i$
(\ref{1PK_def}) to write down the equations of motion:
$$
\rho g_{ab} \left( \frac{D_g \dot{\phi}}{D t} \right)^b
= {{{\mathcal P}_a}^i}_{|i}
- \dfrac{\partial P}{\partial x^k} J
\left(\left(\frac{\partial \phi}{\partial x}\right)^{-1}\right)_a^k,
$$
together with the constraint (\ref{EL_const}).
The above equation can be written in a fully covariant form
$$
\rho \frac{D_g V}{D t} = \operatorname{DIV} {\mathcal P}
- \operatorname{grad} p \circ \phi,
$$
where $V = \dot{\phi}$ is the velocity vector field, ${\mathcal P}$
is the first Piola-Kirchhoff stress tensor, and $p$ is the spatial
pressure.

\section{Symmetries, \ Momentum Maps and Noether's \ Theorem}
\label{Symm_sec}

We already mentioned in the introduction that homogeneous fluid
dynamics has a huge  symmetry, namely the particle relabeling
symmetry, while standard elasticity (usually assumed to be inhomogeneous)
has much smaller symmetry groups, such
as rotations and translations in the Euclidean case.
While inhomogeneous fluids (especially the compressible ones)
are of great interest, the results worked out in \S\ref{Symm_subsec} only
apply to homogeneous fluid dynamics, when the symmetry group is the full
group of volume-preserving diffeomorphisms ${\mathcal D}_\mu$.
However, these results can be generalized to inhomogeneous fluids,
in which case the symmetry group is a \emph{subgroup}
${\mathcal D}^{\rho}_\mu \subset {\mathcal D}_\mu$ that preserves the level
sets of the material density for barotropic fluids, or a \emph{subgroup}
${\mathcal D}^{\rho, \text{ent}}_\mu \subset {\mathcal D}_\mu$ that preserves
the level sets of the material density and entropy for isentropic fluids.

A general model of continuum mechanics will have the
metric $g$ isometry as its symmetry. In particular, the group of rotations
and translations is a symmetry  for models of fluid dynamics and
elasticity in Euclidean spaces. The later is treated in
Marsden, Pekarsky, Shkoller, and West [2000], where the overall
emphasis is on continuum mechanics in Euclidean spaces.

The only symmetry which is universal for \emph{non-relativistic}
continuum mechanics is the time translation invariance. This is due to
the fact that the base manifold is a tensor product of the spatial
part and the time direction, rather than a space-time, so that
all material quantities, such as density $\rho$, metric $G$, etc.
depend only on $x \in B \subset X$.
In this section we shall treat these symmetries separately.
We start with the particle relabeling symmetry, introducing the
necessary notations.

\subsection{Relabeling Symmetry of Ideal Homogeneous Hydrodynamics}
\label{Symm_subsec}

In this subsection we shall consider both the barotropic model
and the incompressible model of ideal homogeneous fluids with fixed boundaries
at the same time.
Their corresponding Lagrangians differ only by the constraint term and both
are equivariant with respect to the action of the group of
volume-preserving diffeomorphisms.

\paragraph{The Group Action.}

The action of the diffeomorphism group ${\mathcal D}_\mu
(B)$ on the (spatial part of the) base manifold $B
\subset X$ captures precisely the meaning of particle
relabeling. For any $\eta \in {\mathcal D}_\mu (B)$, denote this
action by $\eta_X : (x,t) \mapsto (\eta(x), t)$. The lifts of
this  action to the bundles $Y$ and $E$ are given by
$\eta_Y : (x,t,y) \mapsto (\eta(x), t, y)$ and
$\eta_E : (x,t, y, \lambda) \mapsto  (\eta(x), t, y, \lambda)$,
respectively. Both lifts are fiber-preserving and act on the
fibers themselves  by the identity transformation. The
coordinate expressions  have the following form:
\begin{equation}
\label{Diff_action}
\eta_X^0 = \operatorname{Id} \cdot t, \ \eta_X^i = \eta^i (x), \quad
\eta_Y^a = \delta^a_b y^b, \quad
\eta_E^a = (\delta^a_b y^b, \operatorname{Id} \cdot \lambda).
\end{equation}

\paragraph{Jet Prolongations.}

The jet prolongations are natural lifts of automorphisms of $Y$
to automorphisms of its first jet $J^1 Y$ and can be viewed
as covariant analogues of the tangent maps (see Gotay, Isenberg
and Marsden [1997]).

Let $\gamma$ be an element of $J^1 Y$ and $\bar{\gamma}$ be a
corresponding element of the extended phase phase space $J^1 E$,
in coordinates $\gamma = (x^\mu, y^a, v^a_\mu)$ and
$\bar{\gamma} = (x^\mu, y^a, \lambda, v^a_\mu, \beta_\mu)$.
The prolongation of $\eta_Y$ is defined by
\begin{equation}
\label{jet_prol}
\eta_{J^1 Y} (\gamma) = T \eta_Y \circ \gamma \circ T \eta_X^{-1},
\qquad
\eta_{J^1 E} (\bar{\gamma}) = T \eta_E \circ \bar{\gamma} \circ T \eta_X^{-1}.
\end{equation}
We shall henceforth consider $\eta_{J^1 E}$, since it includes
$\eta_{J^1 Y}$ as a special case. In coordinates, we have
$$
\eta_{J^1 E} (\bar{\gamma}) =
\left(\eta^k(x), t; y^b , \lambda ;
{v^a}_0, {v^a}_m \left(
(\frac{\partial \eta}{\partial x})^{-1} \right)_j^m ;
\beta_0, \beta_m \left(
(\frac{\partial \eta}{\partial x})^{-1} \right)_j^m
\right).
$$

If $\xi$ is a vector field on $E$ whose
flow is $\eta_\epsilon$, then its
prolongation $j^1 \xi$ is the vector field on
$J^1 E$ whose flow is $j^1(\eta_\epsilon)$, that is
$j^1 \xi \circ j^1 (\eta_\epsilon)
= (d/d\epsilon) j^1 (\eta_\epsilon)$.
In particular, the vector field $\xi$ corresponding to $\eta_E$ given
by (\ref{Diff_action}) has coordinates $(\xi^i,0,0,0)$ and is
divergence-free; its prolongation $j^1 \xi$, which
corresponds to the prolongation $\eta_{J^1 E}$ of $\eta_E$,
has the following coordinate expression:
\begin{equation}
\label{Diff_action_ext}
j^1 \xi = \left( \xi^i,0; 0,0;0,-{v^a}_m \frac{\partial \xi^m}{\partial x^j};
0, -\beta_m\frac{\partial \xi^m}{\partial x^j} \right).
\end{equation}

\paragraph{Noether's  Theorem.}

Suppose the Lie group ${\mathcal G}$ acts on $\mathcal{C} $ and leaves
the action $S$ invariant.
This is equivalent to the  Lagrangian density $\mathcal{L} $ being
{\bfi equivariant} with respect to ${\mathcal G}$, that is, for all
$\eta \in {\mathcal G}$ and $\gamma \in J^1 Y$,
$$
\mathcal{L}
(\eta_{J^1 Y} (\gamma)) = (\eta_X^{-1})^* \mathcal{L}  (\gamma),
$$
where $(\eta_X^{-1})^* \mathcal{L}  (\gamma)
= (\eta_X)_* \mathcal{L}  (\gamma)$
is a push-forward; this equality means equality of $(n+1)$-forms
at $\eta (x)$.
Denote the {\bfi covariant momentum map} on $J^1 Y$ by
$J_\mathcal{L} \in L(\mathfrak{g}, \Lambda^n(J^1 Y))$.
It is defined by the following expression
\begin{equation}
j^1(\xi) \intprod \Omega_\mathcal{L}  =
d J_\mathcal{L} (\xi)
\label{s22}
\end{equation}
and can be thought of as a Lie algebra valued $n$-form on $J^1 Y$.

Recall that $\phi$ is a solution of the Euler-Lagrange equations
if and only if
$$
(j^1 \phi )^* ({\mathcal W} \intprod \Omega_{{\mathcal L}}) = 0
$$
for any vector field ${\mathcal W}$ on $J^1 Y$. In particular,
setting ${\mathcal W} = j^1 (\xi)$ and applying $(j^1 \phi)^*$
to (\ref{s22}), we obtain the following basic Noether
conservation law:
\begin{thm}
\label{Noether_thm}
Assume that group ${\mathcal G}$ acts on $Y$ by $\pi_{XY}$-bundle
automorphisms and that the Lagrangian density $\mathcal{L}$ is
equivariant with respect to this action for any $\gamma \in J^1 Y$.
Then, for each $\xi \in {\mathfrak g}$
\begin{equation}
\label{Noth_thm}
\mathbf{d} \left( (j^1 \phi)^* J_{\mathcal{L}} (\xi) \right) = 0
\end{equation}
for any section $\phi$ of $\pi_{XY}$ satisfying the Euler-Lagrange equations.
The quantity $(j^1 \phi)^* J_{\mathcal{L}} (\xi)$ is called the
{\bfi Noether current}.
\end{thm}
See Gotay, Isenberg and Marsden [1997]  for a proof.

\paragraph{The Variational Route to Noether's  Theorem.}
The variational route to the covariant Noether's theorem on $J^1 Y$
was first presented on pages 374--375 in
Marsden, Pat\-rick, and Shkoller [1998].  We shall briefly describe this
formulation now.

Recall the notations of the maps $\varphi : U \rightarrow Y$
and the corresponding induced local sections
$\varphi \circ \varphi^{-1}_X$ of $Y$ from Subsection \ref{LDyn_subsec}.
Here again it is important to allow for both vertical and horizontal
variations of the sections.
Vertical variations alone capture only fiber preserving symmetries
(i.e., spatial symmetries), while taking arbitrary variations allows
for both material and spatial symmetries to be included.

The invariance of the action  $S = \int_{U_X} \mathcal{L} $ under the
Lie group action is formally represented by the following expression:
\begin{equation} \label{s22a}
S(\eta_{Y}\cdot \varphi) = S (\varphi) \ \;
\text{for all} \; \
\eta_{Y} \in {\mathcal G}.
\end{equation}
Equation (\ref{s22a}) implies that for each $\eta_{Y}\in {\mathcal G}$,
$\eta_{Y} \cdot \varphi$ is a solution of the Euler-Lagrange equations,
whenever $\varphi$ is a solution.  We
restrict the action of ${\mathcal G}$ to the space of solutions, and let
$\xi_\mathcal{C} $ be the corresponding infinitesimal generator
on ${\mathcal C}$ restricted to the space of solutions; then
\begin{eqnarray*}
0 = (\xi_\mathcal{C}  \intprod d{\mathcal S}) (\varphi) &=&
\int_{\partial U_X} j^1(\varphi \circ \varphi^{-1}_X )^*
[j^1(\xi) \intprod \Theta_\mathcal{L} ] \\
&=& \int_{U_X} j^1(\varphi \circ \varphi^{-1}_X )^*
[j^1(\xi) \intprod \Omega_\mathcal{L} ],
\end{eqnarray*}
since the Lie derivative ${\mathfrak L}_{j^1(\xi)} \Theta_\mathcal{L}  = 0$
by (\ref{s22a}) and Corollary \ref{cor4.2}.

Using (\ref{s22}), we find that
$\int_{U_X} d [ j^1(\varphi \circ \varphi^{-1}_X )^* J_\mathcal{L} (\xi)] =0$,
and since this holds for arbitrary regions $U_X$,
the integrand must also vanish.
Recall that $\phi = \varphi \circ \varphi_X^{-1}$
is a true section of the bundle $Y$, so that this is precisely a
restatement of the Noether's Theorem \ref{Noether_thm}.

\paragraph{Covariant Canonical Transformations.}

The computations of the momentum map from definition
(\ref{s22}) can be simplified significantly in some special cases
which we discuss here. A  $\pi_{X J^1 Y}$-bundle
map $\eta_{J^1 Y} : J^1 Y \rightarrow J^1 Y$ covering the diffeomorphism
$\eta_X : X \rightarrow X$ is called a
{\bfi covariant canonical transformation} if
$\eta_{J^1 Y}^* \Omega_\mathcal{L}  = \Omega_\mathcal{L} $.
It is called a {\bfi special covariant canonical
transformation} if $\eta_{J^1 Y}^* \Theta_\mathcal{L}  = \Theta_\mathcal{L} $.
Recall that forms $\Omega_\mathcal{L} $ and $\Theta_\mathcal{L} $ can
be obtained either by variational arguments or by pulling back canonical
forms $\Omega$ and $\Theta$ from the dual bundle using the Legendre
transformation ${\mathbb F} \mathcal{L} $.

  From Gotay, Isenberg and Marsden [1997], any $\eta_{J^1 Y}$ which
is obtained by lifting some action $\eta_Y$ on $Y$ to $J^1Y$,
is automatically a special canonical transformation.
In this case the momentum mapping is given by
\begin{equation}
\label{mom_map_sp_cov}
J_\mathcal{L} (\xi) =
j^1\xi \intprod \Theta_\mathcal{L} .
\end{equation}

We remark that the validity of this expression does not rely on
the way in which the Cartan form was derived, i.e., for simplicity
of the computations in concrete examples, one can forgo the
issues of vertical vs. arbitrary variations in the variational
derivation and obtain the Cartan form directly from the dual
bundle by means of Legendre transformations. Then, evaluating
this form on the prolongation of a vector of an infinitesimal
generator gives the momentum $n$-form.

\paragraph{Equivariance of the Lagrangian.}
To apply Theorem \ref{Noether_thm} to our case we need to establish
equivariance of the fluid Lagrangians:
\begin{prop}
\label{prop_5.1}
The Lagrangian of an ideal homogeneous barotropic fluid
(\ref{Lag_def})
and the Lagrangian of an ideal homogeneous incompressible fluid
(\ref{Lag_incm_def})
are  equivariant with respect to the ${\mathcal D}_\mu(B)$
action  (\ref{Diff_action}):
$$
\mathcal{L}  ( \eta_{J^1 E} (\bar{\gamma})) =
(\eta_X^{-1})^* \mathcal{L}  (\bar{\gamma}),
$$
for all $\bar{\gamma} \in J^1 E$.
\end{prop}

\begin{proof}
First observe that the material density of an ideal homogeneous
(compressible or incompressible) fluid is constant.
Notice also that Lagrangians (\ref{Lag_def}) and (\ref{Lag_incm_def})
differ only in the potential energy terms. Both these terms are
functions of the Jacobian, which is equivariant with respect to
the action of volume preserving diffeomorphisms given by
(\ref{Diff_action_ext}). Indeed,
\begin{align*}
f \left(J (\eta_{J^1 E} (\bar{\gamma})) \right) d^{n+1} x & =
f \left(\frac{\sqrt{\det [g]}}{\sqrt{\det [G]}} \det (v)
\det \left( \frac{\partial \eta}{\partial x}
\right)^{-1} \right) d^{n+1} x \\
& = (\eta_X^{-1})^* \left( f (J (\bar{\gamma})) d^{n+1} x \right)
\end{align*}
due to the fact that $\det \partial_i \eta^j = 1$ for a volume
preserving diffeomorphism $\eta$; here $f$ can be any function,
e.g. the stored energy $W$ or the constraint $\Phi$.

For the same reason, and the fact that (\ref{Diff_action_ext})
acts trivially on ${v^a}_0$, the kinetic part of both Lagrangians
is also equivariant.
\end{proof}

Proposition \ref{prop_5.1} enables us to use (\ref{mom_map_sp_cov})
for explicit computations of the momentum maps for the relabeling
symmetry of homogeneous hydrodynamics.
We shall consider barotropic and incompressible ideal fluids
separately because their Lagrangians and, hence, their momentum
mappings are different.

\subsubsection*{Barotropic Fluid}

Using (\ref{mom_map_sp_cov}) we can compute the Noether current
corresponding to the relabeling symmetry of the Lagrangian
(\ref{Lag_def}) to be
\begin{multline}
\label{Noet_curr_def}
j^1 (\phi)^* J_{\mathcal{L}} (\xi)  =
\left( \frac{1}{2} \rho  g_{ab} \dot{\phi}^a \dot{\phi}^b
- \rho W - P J \right)  \sqrt{\det [G]} \xi^k d^n x_k - \\
\left( g_{ab} \dot{\phi}^b \phi^a_{,k}\right)
\rho \sqrt{\det [G]} \xi^k d^n x_0,
\end{multline}
where $j^1 \xi$ is the prolongation of the vector field $\xi$
and is given by (\ref{Diff_action_ext}).

The differential of this quantity restricted to the solutions of
the Euler-Lagrange equation is identically zero according to
Theorem \ref{Noether_thm}. Conversely, requiring the differential
of (\ref{Noet_curr_def}) to be zero for arbitrary sections $\phi$ recovers
the Euler-Lagrange equation.
Indeed, computing the exterior derivative and taking into account that
the vector field $\xi$ is divergence free, we obtain:
$$
g_{ab} \left( \frac{D_g \dot{\phi}}{D t} \right)^b = -
\frac{\partial P}{\partial x^k} J
\left(\left(\frac{\partial \phi}{\partial x}\right)^{-1}\right)_a^k,
$$
which coincides with the Euler-Lagrange equation (\ref{EL_bf}).

\subsubsection*{Incompressible Ideal Fluid}

Similar computations using  Lagrangian (\ref{Lag_incm_def})
with the potential energy set to a constant
give the following expression for the Noether current
corresponding to the relabeling symmetry:
\begin{multline}
\label{Noet_curr_incm_def}
j^1 (\phi)^* J_{\mathcal{L}} (\xi)  =
\left( \frac{1}{2} \rho  g_{ab} \dot{\phi}^a \dot{\phi}^b
- P \right)  \sqrt{\det [G]} \xi^k d^n x_k - \\
\left( g_{ab} \dot{\phi}^b \phi^a_{,k}\right)
\rho \sqrt{\det [G]} \xi^k d^n x_0.
\end{multline}

The assumptions of Theorem \ref{Noether_thm} are satisfied; hence
the exterior differential of this Noether current
$d \left( j^1 (\phi)^* J_{\mathcal{L}} (j^1 \xi) \right)$
is equal to zero for all section $\phi$ which are solutions
of the Euler-Lagrange equations.

Now consider the inverse statement. That is,
let us analyze whether the Noether conservation law implies
the Euler-Lagrange equations for incompressible ideal fluids.
Computing the exterior differential of (\ref{Noet_curr_incm_def})
for an arbitrary section $\bar{\phi} = (\phi, \lambda)$
we obtain:
$$
g_{ab} \left( \frac{D_g \dot{\phi}}{D t} \right)^b = -
\frac{\partial P}{\partial x^k}
\left(\left(\frac{\partial \phi}{\partial x}\right)^{-1}\right)_a^k.
$$

Here, we have used the fact that $\xi$ is a divergence free vector field
on $X$. This is precisely the Euler-Lagrange equation
(\ref{EL_if_incm_a}) with the constraint $J=1$ substituted in it.
We point out that the above equation is not equivalent to the
Euler-Lagrange equations, i.e. the constraint cannot be
recovered from the Noether current.
Notice also that the Noether currents (\ref{Noet_curr_def}) and
(\ref{Noet_curr_incm_def}) are different due to the difference in the
corresponding Lagrangians.


\subsection{Time Translation Invariance}

Lagrangian densities (\ref{Lag_def}) and (\ref{Lag_incm_def})
are equivariant with respect to the group $\mathbb{R} $
action on $Y$, given by $\tau_Y  : (x,t,y) \mapsto (x, t+\tau, y)$
for any $\tau \in \mathbb{R} $. This is because the Lagrangians
are explicitly time independent.  One can readily compute the
jet prolongation of the corresponding infinitesimal generator
vector field $\zeta_Y = (0,\zeta, 0)$, where $\tau = \exp \zeta$.
Then, the pull-back by $j^1 \phi$ of the covariant momentum map
corresponding to this symmetry, which we denote by $J^t_{\mathcal{L} }$
to distinguish it from expressions in the previous section,
is given by the following
$n$-form on $X$:
\begin{multline}
\nonumber
(j^1 \phi)^* J^t_{\mathcal{L} } (\zeta) = \zeta \left( L (j^1 \phi) d^n x_0 -
{p_a}^\mu (j^1 \phi) \dot{\phi}^a d^n x_\mu \right) = \\
   - \zeta \left(  e (j^1 \phi) d^n x_0 +
{p_a}^j (j^1 \phi) \dot{\phi}^a d^n x_j \right)(j^1 \phi),
\end{multline}
where, in the last equality,  we have used the definition of the energy
density $e$ given by (\ref{E_def}).

Noether's Theorem \ref{Noether_thm} implies that the exterior
derivative of this expression will be zero along solutions of the
Euler-Lagrange equations. Computing this divergence for an arbitrary
$\zeta$ recovers the energy continuity equation. For a barotropic
fluid, it is given by
$$
\dot{e} = -\sqrt{\det [G]} \operatorname{DIV}
\left( P J \left( \left(\frac{\partial \phi}{\partial x}\right)^{-1}
\right)_a^j
\dot{\phi}^a \right),
$$
while for standard elasticity the equation has the form:
$$
\dot{e} = \sqrt{\det [G]} \operatorname{DIV}
({{\mathcal P}_a}^j \dot{\phi}^a).
$$
The expressions for an incompressible fluid and elastic medium are
similar.

Alternatively, one can consider the inverse statement and
\emph{require} that $d \left( J^t_{\mathcal{L} } (\zeta) \right) = 0$.
This forces the energy continuity equation to be satisfied for some
arbitrary section $\phi$.

\section{Concluding Remarks and Future Directions}
\label{Disc_sec}

In the last section of our paper we would like to comment on
the work in progress and point out general future directions of
the multisymplectic program. Some of the aspects discussed here
are analyzed in detail in our companion paper Marsden, Pekarsky,
Shkoller, and West [2000].

\paragraph{Other Models of Continuum Mechanics.}

The formalism set up in this paper naturally includes other models
of three-dimensional linear and non-linear elasticity and fluid dynamics,
as well as rod and shell models. For elasticity, the choice of
the stored energy $W$ determines a particular model with the
corresponding Euler-Lagrange equation given by (\ref{EL_cm});
this is a PDE to be solved for the deformation field $\phi$.
Introducing the first Piola-Kirchhoff stress tensor ${\mathcal P}$,
the same equation can be written in a compact fully covariant
form (\ref{EL_el}).
An explicit form of the Euler-Lagrange equations and conservation laws
for rod and shell models are not included in this paper but
can be easily derived by following the steps outlined above.
The constrained director models which are common in such models are
handled well by the formulation of constraints that we use in
\S\ref{Const_sec}.

\paragraph{Constrained Multisymplectic Theories.}

The issue of holonomic vs. non-holonomic constraints in classical
mechanics has a long history in the literature.
Though there are still many open questions, the subject of linear and
affine non-holonomic constraints is relatively well-understood
(see, e.g. Bloch, Krishnaprasad, Marsden, and Murray [1996]).
We already mentioned in \S\ref{Const_sec} that this topic
is wide open for multisymplectic field theories, partly due
to the fact that there is simply no well-defined notion of
a non-holonomic constraint for such theories -- it appears
that one needs to distinguish between time and space partial
derivatives.

As all of the examples under present consideration
are non-relativistic and do not have
constraints involving time derivatives, we used the
restriction of Hamilton's principle to the space of allowed
configurations to derive the equations of motion. Note
that this reduces to vakonomic mechanics in the case of an ODE system
with non-holonomic constraints, and is thus incorrect. Of course a
multisymplectic approach to non-holonomic field theories (such as
one elastic body rolling, while deforming, on another, such as a
real automobile tire on pavement), would be of considerable interest
to develop.

\paragraph{Multisymplectic Form Formula and Conservation Laws.}

A very important aspect of any multisymplectic field theory is
the existence of the multisymplectic form formula (\ref{MSFF})
which is the covariant analogue of the fact that the flow of
conservative systems consists of symplectic maps. We deliberately
avoid here any detailed analysis of the implications of  this
formula to the multisymplectic continuum mechanics and refer the
reader to Marsden, Pekarsky, Shkoller, and West [2000], where it is
treated in the context of Euclidean spaces and discretization.
Preliminary results indicate, however, that applications of the
multisymplectic form formula not only can be linked to some known
principles in elasticity (such as the Betti reciprocity principle),
but also can produce some new interesting relations which depend on the
space-time direction of the first variations $\mathcal{V} ,
\mathcal{W} $ in (\ref{MSFF}). An accurate and consistent
discretization of the model then  results in so called {\bfi
multisymplectic integrators} which preserve  the discrete analogues
of the multisymplectic form and the conservation laws.

\paragraph{Discretization.}

This is another very interesting and important part of our project
which is addressed in detail in our companion paper Marsden, Pekarsky,
Shkoller, and West [2000], where the approach of finite elements
for models in Euclidean spaces
is adopted. It is shown that the finite element method for static
elasticity is a multisymplectic integrator.
Moreover, based on the result in Marsden and West [2000],
it is shown that the finite elements time-stepping with
the Newmark algorithm is a multisymplectic discretization.

As we mentioned in the previous paragraph, a consistent discretization
based on the variational principle would preserve the discrete
multisymplectic form formula together with the discrete multi-momentum maps
corresponding to the symmetries of a particular system.
Then, the integral form of the discrete Noether's
Theorem implies that a sum of the values of the discrete
momentum map over some set of nodes is zero. One implication of
this statement for incompressible fluid dynamics is a discrete
version of the vorticity preservation. Such discrete conservations
are among the hot topics of the ongoing research.

\paragraph{Symmetry Reduction.}

In the previous section we discussed at length the particle
relabeling symmetry of ideal homogeneous hydrodynamics and its multisymplectic
realization. Reduction by this symmetry takes us from the Lagrangian
description in terms of \emph{material} positions and velocities to the
Eulerian description in terms of \emph{spatial} velocities. In the
compressible case one only reduces by the subgroup of the particle
relabeling group that leaves the stored energy function invariant;
for example, if the stored energy function depends on the
deformation only through the density and entropy, then this means
that one introduces them as dynamic fields in the reduction
process, as in Euler-Poincar\'e theory (see Holm, Marsden and Ratiu
[1998].)

In the unconstrained (i.e., defined on the extended jet bundle $J^1 E$)
multisymplectic description of ideal incompressible fluids,
the multisymplectic reduced space is realized as a fiber bundle
$\Upsilon$ over $X$ whose fiber coordinates include the Eulerian
velocity $u$  and some extra field corresponding to compressibility.
Then, the reduced  Lagrangian density determines, by means of
a constrained variational principle,
the Euler-Lagrange equations which give the evolution of the
spatial velocity field $u(x) \in \Upsilon_x$ together with a
condition of $u$ being divergence-free. A general
Euler-Poincar\'{e} type theorem relates this equation with equation
(\ref{EL_if_incm2}) by relating the corresponding variational
principles.

Such a description is a particular example of a general  procedure
of multisymplectic reduction. The case of a finite-dimensional
vertical group action was first considered
in Castrill\'{o}n-L\'{o}pez, Ratiu and Shkoller [2000]. More
general cases of an infinite-dimensional group action such as that
for incompressible ideal hydrodynamics, electro-magnetic fields and
symmetries in complex fluids is planned for a future publication.
The reader is also referred to a related work by  Fern\'{a}ndez,
Garc\'{i}a, and Rodrigo [1999].

\paragraph{Vortex Methods.}

One of our ultimate objectives is to further develop, using
the multisymplectic approach, some methods and techniques
which were derived in the infinite-dimensional framework
and which proved to be very useful. One of them is the
vortex blob method developed by Chorin [1973], which
recently has been linked to the so-called averaged Euler
equations of ideal fluid (see Oliver and Shkoller [1999]).

\paragraph{Higher Order Theories.}
Constraints involving higher than first-order derivatives are beyond
the current exposition and should be treated in the context of
higher-order multisymplectic field theories defined on $J^k Y, \ k
>   1$.

The averaged Euler equations (see Holm, Marsden and Ratiu [1998] and
Marsden, Ratiu and Shkoller [2000] and references therein) provide
an interesting example of a higher order fluid theory with
constraints (depending only on first derivatives of the field) to
which the multisymplectic methods can presumably be applied by
using the techniques of Kouranbaeva and Shkoller [2000]. It would be
interesting to carry this out in detail. In the long run, this
promises to be an important computational model, so that its
formulation as a multisymplectic field theory and the
multisymplectic discretization of this theory is of considerable
interest.

\paragraph{Covariant Hamiltonian description.}

Finally, another very interesting aspect of the project is
developing the multi-Hamiltonian description of continuum mechanics
along the lines outlined in Marsden and Shkoller [1999].

\subsection*{Acknowledgments}
The authors would like to thank Sanjay Lall, Tudor Ratiu and Michael
Ortiz for their comments and interest in this work.

\subsection*{References}

\begin{description}

\item Arnold, V.I. [1966]
Sur la g\'{e}om\'{e}trie differentielle
des groupes de Lie de dimenson
infinie et ses applications \`{a}
l'hydrodynamique des fluids parfaits.
{\it Ann. Inst. Fourier, Grenoble\/} {\bf 16}, 319--361.

\item Arnold, V.I. and B. Khesin [1998]
{\it Topological Methods in Hydrodynamics.}
Appl. Math. Sciences {\bf 125}, Springer-Verlag.

\item Bloch, A. M. and P.E. Crouch, [1999] Optimal control,
optimization, and analytical mechanics. {\it Mathematical control
theory}, 268--321, Springer, New York, 1999.

\item Bloch A.M. and P.S. Krishnaprasad and J.E. Marsden and R.M.
Murray [1996] Nonholonomic mechanical systems with symmetry,
{\it Arch. Rat. Mech. An.}, {\bf 136}, 21--99

\item  Castrill\'{o}n-L\'{o}pez, M., T.S. Ratiu and  S. Shkoller [2000]
Reduction in Principal Fiber Bundles:
Covariant Euler-Poincar\'{e} Equations, {\it Proc. Am. Math.
Soc.}, to appear.

\item Chorin, A. [1973]
Numerical study of slightly viscous flow,
{\it J. Fluid Mech.} {\bf 57}, 785-796.

\item Ebin, D.G. and J.E. Marsden [1970]
Groups of diffeomorphisms and the motion of an incompressible
fluid, {\it Ann. Math.\/} {\bf 92}, 102--163.

\item Fern\'{a}ndez A., P.L. Garc\'{i}a, and C. Rodrigo [1999]
Stress-energy-momentum tensors in higher order variational calculus,
{\it J. Geom. Phys.}, to appear.

\item Gotay, M., J. Isenberg, and J.E. Marsden [1997]
{\it Momentum Maps and the Hamiltonian Structure of
Classical Relativistic Field
Theories, I and II,\/} available from
\textcolor{blue}{\url{http://www.cds.caltech.edu/~marsden/}}.

\item Gotay, M.J. and J.E. Marsden [1992]
Stress-energy-momentum tensors and the
Belifante-Resenfeld formula.
{\it Cont. Math. AMS.\/} {\bf 132}, 367--392.

\item Holm, D.D., J. E. Marsden and T. S. Ratiu [1998]
The Euler--Poincar\'{e} equations and semidirect products
with applications to continuum theories,
{\it Adv. in Math.}, {\bf 137}, 1-81.

\item Jalnapurkar S.M. [1994] Modeling of Constrained Systems, \\
\textcolor{blue}{\url{http://www.cds.caltech.edu/\~{}smj/}}.

\item Kijowski, J. and Magil, G [1992]
Relativistic elastomechanics as a Lagrangian field theory.
{\it J. Geo. Phys.\/} {\bf 9}, 207-223.

\item Kijowski, J. and W. Tulczyjew [1979]
{\it A Symplectic Framework for Field
Theories.\/} Springer Lect. Notes in Physics {\bf 107}.

\item Kouranbaeva S. and S. Shkoller [2000]
A variational approach to second-order multisymplectic field
theory, {\it J. Geom. Phys.}, to appear.

\item Lewis A.D. [1996]
The geometry of the {G}ibbs-{A}ppell equations and {G}auss's
Principle of Least Constraint,
{\it Rep. on Math. Phys.}, {\bf 38}, 11-28.


\item Lewis A.D. and R.M. Murray [1995]
Variational principles for constrained systems: theory and experiment,
{\it The Intern. J. of Nonlinear Mech.}, {\bf 30},
793--815.

\item Lu J. and P. Papadopoulos [1999]
A covariant constitutive approach to finite plasticity,
to appear in {\it Proc. 8th International Symposium on Plasticity}.

\item D. G. Luenberger [1969] Optimization by Vector Space Methods.
John Wiley, New York.

\item Marsden, J.E. [1988],
The Hamiltonian formulation of classical field theory,
{\it Cont. Math. AMS} {\bf 71}, 221-235.

\item Marsden, J.E. [1992],
{\it Lectures on Mechanics\/} London 
Mathematical
Society Lecture note series, {\bf 174}, Cambridge 
University Press.

\item Marsden, J.E. and T.J.R. Hughes [1983]
{\it 
Mathematical Foundations of Elasticity.\/}
Prentice Hall, reprinted 
by Dover Publications, N.Y., 1994.

\item Marsden, J. E., G. W. 
Patrick, and S. Shkoller [1998]
Multisymplectic Geometry, Variational 
Integrators, and
Nonlinear PDEs, {\it Comm. Math. Phys.} {\bf 199}, 
351--395.

\item Marsden J. E., S. Pekarsky, S. Shkoller and M. West 
[2000]
Multisymplectic continuum mechanics in Euclidean spaces,
{\it 
preprint}.

\item Marsden, J.E., T.S. Ratiu, and S. Shkoller 
[2000]
The geometry and analysis of the averaged Euler equations
and 
a new diffeomorphism group, {\it Geom. and Funct. An.}
(to 
appear).

\item Marsden, J.E. and T.S. Ratiu [1999]  {\it 
Introduction to
Mechanics and Symmetry.\/} Texts in Applied 
Mathematics, {\bf  17},
Springer-Verlag, 1994. Second Edition, 
1999.

\item Marsden J. E. and S. Shkoller [1999]
Multisymplectic 
geometry, covariant Hamiltonians, and water
waves
{\it Math. Proc. 
Camb. Phil. Soc.}, {\bf 125},  553--575.

\item Oliver M. and S. 
Shkoller [1999]
The vortex blob method as a second-grade 
non-Newtonian fluid, E-print, \\
\textcolor{blue}{\url{http://xyz.lanl.gov/abs/math.AP/9910088/}}.

\item Shkoller, S. [1998] Geometry and curvature of diffeomorphism
groups with $H^1$ metric and mean hydrodynamics, {\it
J. Func. Anal.} {\bf  160}, 337--365.

\end{description}

\end{document}